\documentclass{amsart}

\usepackage{amsmath,amsfonts,amsthm,amssymb,eucal}

\theoremstyle{plain}
\newtheorem{Theorem}{Theorem}[section]
\newtheorem{Corollary}[Theorem]{Corollary}
\newtheorem{Proposition}[Theorem]{Proposition}
\newtheorem{Lemma}[Theorem]{Lemma}
\newtheorem{Sublemma}{Sublemma}

\theoremstyle{definition}
\newtheorem{Remark}{Remark}

\title{Small generators for $S$-unit groups of division algebras}

\author{Ted Chinburg}
\address{Department of Mathematics \\ University of Pennsylvania \\ 209 South 33rd Street \\ Philadelphia, PA 10104}
\email{ted@math.upenn.edu}
\thanks{Chinburg partially supported by NSF Grant DMS 1100355.}

\author{Matthew Stover}
\address{Department of Mathematics \\ Temple University \\ 1805 N.\ Broad Street \\ Philadelphia, PA 19122}
\email{mstover@temple.edu}
\thanks{Stover partially supported by NSF RTG grant DMS 0602191 and NSF grant DMS 1361000.}



\begin{document}

\maketitle

\begin{center}
To H.\ W.\ Lenstra, Jr.
\end{center}

%%%%%%%%%%%%%%%%%%%%

\begin{abstract}
Let $k$ be a number field, suppose that $B$ is a central simple division algebra over $k$, and choose any maximal order $\mathcal{D}$ of $B$. The object of this paper is to show that the group $\mathcal{D}_S^*$ of $S$-units of $B$ is generated by elements of small height once $S$ contains an explicit finite set of places of $k$. This generalizes a theorem of H.\ W.\ Lenstra, Jr., who proved such a result when $B = k$. Our height bound is an explicit function of the number field and the discriminant of a maximal order in $B$ used to define its $S$-units.
\end{abstract}

%%%%%%%%%%%%%%%%%%%%

\section{Introduction}\label{intro}

%%%%%%%%%%%%%%%%%%%%

Let $k$ be a number field, suppose that $B$ is a central simple division algebra over $k$, and choose any maximal order $\mathcal{D}$ of $B$. The object of this paper is to show that the group $\mathcal{D}_S^*$ of $S$-units of $k$ is generated by elements of small height once $S$ contains an explicit finite set of places of $k$. 

A result of this kind was shown by Lenstra in \cite{Lenstra} when $B$ is $k$ itself. In this case, $G$ is the multiplicative group $\mathbb{G}_m$ and  the notion of height is the classical one. Lenstra showed that once $S$ is sufficiently large, the group $\mathbb{G}_m(O_{k,S}) = O_{k,S}^*$ is generated by elements whose log heights are bounded by
\[
\frac{1}{2} \log |d_{k/\mathbb{Q}}| + \log m_S + r_2(k) \log(2/\pi),
\]
where $d_k$ is the discriminant of $k$, $m_S$ is the maximal norm of a finite place in $S$, and $r_2(k)$ is the number of complex places of $k$.  One version of our results is the following.

%%%%%%%%%%%%%%%%%%%%

\begin{Theorem}\label{bounded height intro}
Suppose $B$ is a central simple division algebra of dimension $d^2$ over a number field $k$, $n = [k : \mathbb{Q}]$, and $s$ is the number of real places of $k$ over which $B$ ramifies. Then there is a maximal order $\mathcal{D}$ of $B$ with discriminant $d_{\mathcal{D}}$ and functions $f_1(n, d)$ and $f_2(n, d)$ of integer variables $n$ and $d$ for which the following is true. Define
\[
e = \frac{2n}{d(2n-s)}.
\]
Suppose that $S$ is a finite set of places of $k$ containing all the archimedean places and that $S$ contains all finite places $v$ such that
\[
\mathrm{Norm}(v) \le f_1(n, d)\, d_{\mathcal{D}}^e.
\]
Let $m_{S_f}$ be the maximum norm of a finite place in $S$. Then $e \le 1$ and the group $\Gamma_S$ of $S$-units in $B$ with respect to the order $\mathcal{D}$ is generated by the finite set of elements of height bounded above by
\[
f_2(n, d)\, m_{S_f}\, d_{\mathcal{D}}^e.
\]
\end{Theorem}

%%%%%%%%%%%%%%%%%%%%

See the remarks at the end of \S \ref{ssec:explicit} for explicit expressions for $f_1(n, d)$ and $f_2(n, d)$. In many cases (e.g., when $B$ is a number field) we have
\begin{equation}\label{eq:introf1def}
f_1(n, d) = d^{n d} \left( \frac{2}{\pi} \right)^{\frac{n d r_2}{n - s / 2}} \left( \frac{2 \sqrt{2}}{\pi} \right)^{\frac{n d s}{2 n - s}}
\end{equation}
\begin{equation}\label{eq:introf2def}
f_2(n, d) = d^{n d + n + s} \left( (d - 1)! \right)^{n - s} 2^{s \frac{d^2 - 2 d - 4}{2}} \left( \frac{2}{\pi} \right)^{\frac{n d r_2}{n - s / 2}} \left( \frac{2 \sqrt{2}}{\pi} \right)^{\frac{n d s}{2 n - s}}
\end{equation}
where $r_2$ is the number of complex places of $k$ and all other notation is from the statement of Theorem \ref{bounded height intro}. Thus, when $B$ is a number field this exactly reproduces Lenstra's bound.

When $S$ is empty, even for a number field $k$ one does not expect to be able to generate the unit group $O_k^*$ by elements whose log heights are bounded by a polynomial in $\log|d_{k/\mathbb{Q}}|$. For example, the Brauer--Siegel Theorem implies that if $k$ is real quadratic of class number $1$, then the log height of a generator of $O_k^*$ is greater than $c_\epsilon \, d_{k/\mathbb{Q}}^{1/2 - \epsilon}$ for all $\epsilon > 0$, where $c_\epsilon > 0$ depends only on $\epsilon$. However, to our knowledge there is no unconditional proof, even in the case of real quadratic fields, that there cannot be an upper bound on the log heights of generators for $O_k^*$ that is polynomial in $\log |d_{k/\mathbb{Q}}|$.

To develop our counterpart of Lenstra's results, we must first define an intrinsic notion of height for elements of $B^*$. The role of $|d_{k / \mathbb{Q}}|$ is played by the discriminant $d_{\mathcal{D}}$ of an $O_k$-order $\mathcal{D}$ in $B$ that is used to define the $S$-integrality of points of $G$ over $k$. The height bound we produce applies to all choices of $\mathcal{D}$ once $S$ is sufficiently large. It implies, in particular, that there is a maximal order $\mathcal{D}$ in $B$ such that when $S$ is sufficiently large (in an explicitly defined sense), one can generate the points of $G$ over $O_{k,S}$ by elements whose log heights are bounded by $\frac{1}{2} \log |\Delta_{k/\mathbb{Q}}| + \log m_S + \mu$ where $\mu$ depends only only the degree of $B$ over $\mathbb{Q}$. See \cite{Liebendorfer}, \cite{Liebendorfer2}, \cite{ChanFukshansky} and references therein for work on heights on quaternion algebras.

To show that our bounds remain quite effective outside the number field setting, in \S \ref{sec:Hamilton} we apply our results to Hamilton's quaternion algebra over $\mathbb{Q}$. The main result is the following.

%%%%%%%%%%%%%%%%%%%%

\begin{Theorem}\label{thm:HamiltonIntro}
Let $B$ be Hamilton's quaternion algebra over $\mathbb{Q}$, that is, the rational quaternion algebra with basis $\{1, I, J, IJ\}$ such that $I^2 = J^2 = -1$ and $I J = - J I$. Let $\mathcal{D}$ be the maximal order
\[
\mathbb{Z} \left[ 1, I, J, \frac{1 + I + J + I J}{2} \right]
\]
and $S = \{\infty, \ell_1, \dots, \ell_h\}$ be a set of places containing the archimedean place $\infty$ and any set $\{\ell_i\}_{i = 1}^h$ of distinct odd primes. Then the unit group $\mathcal{D}_S^*$ is generated by the finite set of elements with reduced norm contained in $\{1, \ell_1, \dots, \ell_h\}$.
\end{Theorem}

%%%%%%%%%%%%%%%%%%%%

As in Lenstra's case, Theorem \ref{bounded height intro} leads to an algorithm for finding generators for $O_{k,S}^*$. After embedding $B$ into a real vector space, the algorithm is reduced to the classical problem of enumerating lattice points of bounded norm. That an algorithm exists to generate $S$-integral points of an algebraic group over a number field, with no assumptions on $S$, was known by work of Grunewald--Segal \cite{Grunewald--Segal}, \cite{Grunewald--Segal2}. Unlike their algorithm, ours is primitive recursive, which answers a question raised in \cite{Grunewald--Segal2}.

Lenstra went on to show that his algorithm, together with linear algebra, can be used to give a deterministic algorithm for finding generators for the unit group $O_k^*$ with running-time  $O_{\epsilon}(|d_{k / \mathbb{Q}}|^{3/4 +\epsilon})$. The output of the algorithm consists of the digits of a set of generators. Answering a question raised by Lenstra, Schoof recently announced an improvement of this run-time to $O_{\epsilon}(|d_{k / \mathbb{Q}}|^{1/2 +\epsilon})$. By the discussion of the Brauer--Siegel theorem above, one expects the length of the output may be on the order of $|d_{k / \mathbb{Q}}|^{1/2-\epsilon}$ in some cases, but the input to the problem, namely enough information to specify $k$, will in general be much shorter. For instance, a real quadratic field $k$ can be specified by giving its discriminant $d_{k / \mathbb{Q}}$, and the number of bits necessary to specify $d_{k / \mathbb{Q}}$ is proportional to $\log|d_{k / \mathbb{Q}}|$.

Lenstra's algorithm for finding generators for $O_k^*$ makes essential use of the fact that $O_{k,S}^*$ is abelian, and this is not the case when $B$ is noncommutative and $\mathcal{D}_S^*$ is infinite. Consequently, we do not know a counterpart of this algorithm for producing generators in the general case.

We also note that there is a spectral approach to finding small generators for groups acting on symmetric spaces. For example, if $V$ is a compact quotient of a rank one symmetric space other than hyperbolic $2$- or $3$-space, Burger and Schroeder \cite{Burger--Schroeder} showed that one can bound the diameter of $V$ from above in terms of its volume and $\lambda_1(V) \mathrm{diam}(V)$ in terms of $\log \mathrm{vol}(V)$. One can use this to bound the length (in the Riemannian metric on $V$) of a generating set. As noted in \cite{Burger--Schroeder}, these results fail for hyperbolic $2$- and $3$-space, though Peter Shalen informed us that one could prove an analogous theorem for arithmetic Fuchsian and Kleinian groups assuming Lehmer's conjecture. However, the problem we solve in this paper is of a different kind, even where a result like that of \cite{Burger--Schroeder} holds. Our methods find explicit matrix generators with small entries, and it is not at all clear that generators that are short in an associated Riemannian metric have representatives in $\mathrm{GL}_n$ with small entries. It would be interesting to see if spectral methods for $S$-arithmetic subgroups of reductive groups could produce generators that have small height.

The main tool for producing small generators for $O_{k,S}^*$ once $S$ is sufficiently large is Minkowski's lattice point theorem. This determines elements of $B^*$ which can be shown to be $S$-integral by careful consideration of the constants required to apply Minkowski's theorem. In particular, the assumption that $B$ is a division algebra is crucial, and it would be interesting to extend our results to unit groups of arbitrary central simple algebras. The reason we cannot work with a general central simple algebra is because Minkowski's lattice point theorem returns a nonzero element of the algebra $B$. We prove that it is $S$-integral. However, it might not be invertible in $B$ if $B$ is not a division algebra. Given that $\mathrm{SL}_n(\mathbb{Z})$ is generated by elementary matrices, which certainly have small height, we expect such a result to hold. It appears to us that it is a deep problem to extend such height results to generators for the $S$-integral points of more general linear algebraic groups over number fields.

%%%%%%%%%%%%%%%%%%%%

\section{Notation and definitions}\label{notation}

%%%%%%%%%%%%%%%%%%%%

Let $k$ be a number field with ring of integers $O_k$ and $[k : \mathbb{Q}] = n$. We denote by $V_\infty$ (resp.\ $V_f$) the set of archimedean (resp.\ finite) places of $k$. Let $B$ be a division algebra over $k$ with degree $d$ and let $\mathcal{D} \subset B$ be an $O_k$-order in $B$. The multiplicative group of units in a ring $R$ will be denoted $R^*$. For each place $v$ of $k$ and any $k$-algebra or $O_k$-module $A$, let $A_v$ denote the associated completion at $v$.

Define a norm on $B_v^*$ by
\begin{equation}\label{eq:norm-def}
\mathrm{Norm}_v(x_v) = \mathrm{Norm}_{k_v / \mathbb{Q}_{p(v)}}(\det(x_v \curvearrowright B_v)),
\end{equation}
where $p(v)$ is the place of $\mathbb{Q}$ under $v$ and $x_v \curvearrowright B_v$ is the $k_v$-linear endomorphism of $B_v$ induced by left $x_v$-multiplication. If $\mathrm{det}_v:B_v \to k_v$ is the reduced norm, then
\begin{equation}
\label{eq:rednorm}
\det(x_v \curvearrowright B_v) = \mathrm{det}_v(x_v)^d.
\end{equation}

The idele group $J(B)$ is the restricted direct product $\prod^\prime_v B_v^*$ of the $B_v^*$ with respect to the groups $\mathcal{D}^*_v$. For $x = \prod_v x_v \in J(B)$, define
\begin{equation}\label{eq:norm-def-inf}
\mathrm{Norm}_\infty(x) = \prod_{v \in V_\infty} \mathrm{Norm}_v(x_v)
\end{equation}
\begin{equation}\label{eq:norm-def-fin}
\mathrm{Norm}_f(x) = \prod_{v \in V_f} \mathrm{Norm}_v(x_v).
\end{equation}
We view these norms as elements of the idele group $J(\mathbb{Q})$ of $\mathbb{Q}$ in the natural way. Let $|\ |:J(\mathbb{Q}) \to \mathbb{R}_{> 0}$ be the usual norm. By the product formula,
\begin{equation}\label{eq:prod-form}
|\mathrm{Norm}_\infty(x)| = |\mathrm{Norm}_f(x)|^{-1}
\end{equation}
for all $x \in B^*$.

Let $S$ be a finite set of places of $k$ containing $V_\infty$. Set $S_f = S \smallsetminus V_\infty$, and consider the groups
\[
B^*_{\mathbb{R}} = \prod_{v \in V_\infty} B_v^*
\]
\[
B^*_{S_f} = \prod_{v \in S_f} B_v^* \ \subset \ B^*_f = \left. \prod_{v \in V_f} \right.^\prime B_v^*
\]
and the product
\begin{equation}\label{eq:idele-group}
B_S^* = \prod_{v \in S} B^*_v = B^*_{\mathbb{R}} \times B^*_{S_f} \subset J(B) = B^*_{\mathbb{R}} \times B^*_f.
\end{equation}
Let $G_S$ be the subgroup of $B_S^*$ that satisfies the product formula, so
\begin{equation}\label{eq:S-group}
G_S = \left\{ (x, \beta) \in B_S^*\ :\ |\mathrm{Norm}_\infty(x)| = |\mathrm{Norm}_f(\beta)|^{-1} \right\}.
\end{equation}
If $O_{k, S}$ denotes the $S$-integers of $k$, i.e., those elements of $k$ which lie in $O_{k, v}$ for all $v \notin S$, then the $S$-order of $B$ associated with $\mathcal{D}$ is
\begin{equation}\label{eq:S-order}
\mathcal{D}_S = O_{k, S} \otimes_{O_k} \mathcal{D}.
\end{equation}
The group of invertible elements of $\mathcal{D}_S$ will be denoted $\Gamma_S$ and is called the group of $S$-\emph{units} of $\mathcal{D}$. 

We define a topology on $G_S$ by its natural embedding into $B_S^*$. The image of $\Gamma_S$ in $G_S$ under the diagonal embedding is a discrete subgroup. We have diagonal embeddings $\Gamma_S \to B_S^*$ and $\Gamma_S \to \prod_{v \not \in S} \mathcal{D}^*_v$, and the product of these embeddings is the natural diagonal embedding of $\Gamma_S$ into $J(B)$.

For any element
\[
\alpha = \prod_{v \in V_f} \alpha_v \in B_f^*,
\]
there is a right-$\mathcal D$-module
\begin{equation}\label{eq:right-module-def}
\alpha \mathcal{D} = B \cap \left( \prod_{v \in V_f} \alpha_v \mathcal{D}_v \right),
\end{equation}
where $B$ is diagonally embedded in $B_f$. For $\alpha \in \mathcal{D}$, the index of $\alpha \mathcal{D}$ in $\mathcal{D}$ equals $|\mathrm{Norm}_f(\alpha)|^{-1}$. We also have the left-$\mathcal{D}$-module
\begin{equation}\label{eq:left-module-def}
\mathcal{D} \alpha^{-1} = \{ x \in B\ :\ x (\alpha \mathcal{D}) \subseteq \mathcal{D} \}.
\end{equation}

%%%%%%%%%%%%%%%%%%%%

\section{Absolute values and heights}\label{sec:heights}

%%%%%%%%%%%%%%%%%%%%

In this section we define absolute values on the completions $B_v$ of $B$. These will be used to define our notion of height for elements of $B^*$. We point the reader to \cite{Liebendorfer}, \cite{Liebendorfer2}, \cite{ChanFukshansky} and references therein for earlier work on heights for quaternion algebras over number fields.

For each place $v$ of $k$ there is a division algebra $A_v$ over $k_v$ such that $B_v = k_v \otimes_k B$ is isomorphic to a matrix algebra $\mathrm{M}_{m(v)}(A_v)$. The dimension of $A_v$ over $k_v$ is $d(v)^2$ for some integer $d(v)$ such that $d(v) m(v) = d = \sqrt{\dim_k B}$. Note that $A_v$ and $B_v$ have center isomorphic to $k_v$.

For finite $v$ let $O_v$ be the ring of integers of $k_v$. We fix isomorphisms $\rho_v : B_v \to \mathrm{M}_{m(v)}(A_v)$ such that for almost all finite $v$, $A_v = k_v$ and $\rho_v(\mathcal{D}_v) = \mathrm{M}_{m(v)}(O_v)$. Let $\mathrm{N}_v : A_v \to k_v$ be the reduced norm. Then $\mathrm{N}_v(r) = r^{d(v)}$ for any $r \in k_v \subseteq A_v$.

For all places $v$ of $k$, let $|\ |_v$ be the usual normalized absolute value on $k_v$. We extend $|\ |_v$ to an absolute value on $A_v$ by $|\alpha|_v = |\mathrm{N}_v(\alpha)|^{1/d(v)}_v$ for $\alpha \in A_v$. This absolute value is clearly multiplicative and restricts to the usual absolute value on the center $k_v$ of $A_v$.

Suppose $v$ is nonarchimedean. There is a unique maximal order $U_v$ in $A_v$, namely the set of $\alpha \in A_v$ such that $|\alpha|_v \leq 1$. When $A_v \neq k_v$, $U_v$ is a noncommutative local ring, and it is $O_v$ when $A_v = k_v$. The unique maximal two-sided ideal of $U_v$ is the set $P_v$ of $\alpha \in A_v$ for which $|\alpha|_v < 1$. There is an element $\lambda_v$ of $P_v$ such that $P_v = U_v \lambda_v = \lambda_v U_v$; such $\lambda_v$ are called prime elements by Weil in \cite[Def.\ 3, Chap.\ I.4]{Weil}. By \cite[Prop.\ 5, Chap.\ I.4]{Weil}, $\mathrm{N}_v(\lambda_v)$ is a uniformizer in $k_v$. Thus $|\lambda_v|_v = |\mathrm{N}_v(\lambda_v)|_v^{1 / d(v)} = (\# k(v))^{-1 / d(v)}$ where $k(v)$ is the residue field of $v$ and so the range of $|\ |_v$ on $A_v^*$ is $(\# k(v))^{1 / d(v)})^{{\mathbb{Z}}}$. The set of $\alpha \in A_v^*$ such that $|\alpha|_v \le (\# k(v))^{-t / d(v)}$ is exactly $P_v^t$. This implies that $|\alpha + \beta|_v \le \max(|\alpha|_v, |\beta|_v) $ for every $\alpha, \beta \in A_v$.

We now prove a simple lemma.

%%%%%%%%%%%%%%%%%%%%

\begin{Lemma}\label{lem:triangle-cauchy}
With notation as above, suppose that $v$ is archimedean. For all $m$-element subsets $\{\alpha_i\}_{i = 1}^m \subset A_v$,
\begin{equation}
\label{eq:sum}
\left| \sum_{i = 1}^m \alpha_i \right|_v \le m^{[k_v : \mathbb{R}] - 1} \sum_{i = 1}^m |\alpha_i|_v.
\end{equation}
\end{Lemma}

%%%%%%%%%%%%%%%%%%%%

\begin{proof}
If $v$ is complex, then $A_v = k_v \cong \mathbb{C}$ and $| \ |_v$ is the square of the usual Euclidean absolute value. The result reduces to the Cauchy--Schwarz inequality. If $v$ is real and $A_v = k_v \cong \mathbb{R}$, the lemma is again clear.

The final case is when $v$ is real and $A_v$ is Hamilton's quaternions $\mathbb{H} = \mathbb{R} + \mathbb{R}I + \mathbb{R} J + \mathbb{R} IJ$ where $I^2 = J^2 = -1$ and $I J = -J I$. Here,
\[
|a + bI + cJ + dIJ|_v = (a^2 + b^2 + c^2 + d^2)^{1/2}.
\]
We can view this as the Euclidean length in $\mathbb{R}^4$ of the vector $(a, b, c, d)$. It is clear from the triangle inequality that the optimal constant in this case is again $1$. This proves the lemma.
\end{proof}

%%%%%%%%%%%%%%%%%%%%

Recall that for each place $v$ of $k$ we fixed an isomorphism $\rho_v : B_v \to \mathrm{M}_{m(v)}(A_v)$. For each $v$ and each $\gamma \in B_v^* $, let $\gamma^{i,j}(v)$ denote the $(i,j)$-component of the $m(v) \times m(v)$ matrix $\rho_v(\gamma) \in \mathrm{GL}_{m(v)}(A_v)$. Define $|\gamma|_v = \max_{i,j} |\gamma^{i,j}(v)|_v$. Embed $B$ into $B_v = k_v \otimes_k B$ in the natural way. Then the \emph{height} of $\gamma \in B^*$ is defined by
\begin{equation}\label{eq:height-def}
\mathrm{H}(\gamma) = \prod_{v \in V} \max\{1, |\gamma|^{d(v)}_v\},
\end{equation}
where the product is over the set $V = V_\infty \cup V_f$ of all places of $k$. From the definition of $|\ |_v$ we see that 
\begin{equation}
\label{eq:height-def2}
\mathrm{H}(\gamma) = \prod_{v \in V} \max \left\{ 1,\ \max_{i,j} |\mathrm{N}_v(\gamma^{i,j}(v))|_v \right\}.
\end{equation}

For any $S$ and any positive real number $x$, the set
\begin{equation}\label{eq:bdd-ht-def}
\mathrm{BH}_{S}(x) = \{ \gamma \in \Gamma_S\ :\ \mathrm{H}(\gamma) \leq x \}
\end{equation}
of $S$-units of $\mathcal{D}$ with height bounded by $x$ is finite. Indeed, bounding the nonarchimedean height bounds the denominator of each matrix entry under the image of every $\rho_v$, so the set of $\gamma \in \Gamma_S$ with bounded height is contained in a lattice in $B_{\mathbb{R}}$. Bounding the archimedean height immediately implies finiteness.

We end this section by proving some inequalities we will need later concerning the behavior of absolute values on taking products and inverses of elements of $B_v$.

%If $v$ is archimedean, define $\mathrm{det}'_v : B_v \to k_v$ by $\mathrm{det}'_v(q) = |\mathrm{det}_v(q)|_v^{1 / d(v)}$. If $v$ is nonarchimedean, fix a uniformizer $\pi_v$ of $k_v$. Define $\mathrm{det}'_v : B_v \to k_v$ by setting $\mathrm{det}'_v(q) = 0$ if $\mathrm{det}_v(q) = 0$, and by letting $\mathrm{det}'_v(q)$ be the power of $\pi_v$ such that $\mathrm{det}'_v(q)^{d(v)}$ generates the same fractional $O_v$ ideal as $\mathrm{det}_v(q)$. This is well-defined because the $v$-adic valuation of every element of $\mathrm{det}_v(B_v^*)$ is a multiple of $d(v)$.

%%%%%%%%%%%%%%%%%%%%

\begin{Lemma}\label{lem:detbound1}
Suppose $v$ is a finite place of $k$ and fix an isomorphism $\rho_v : B_v \to \mathrm{M}_{m(v)}(A_v)$. Let $\mathrm{det}_v : B_v \to k_v$ be the reduced norm. For $y, y' \in B_v$:
\begin{enumerate}
\label{eq:upperone}
\item[1.] $|y y'|_v \le |y|_v \, |y'|_v$.
\item[2.] If $y$ is invertible, then $|\mathrm{det}_v(y)|_v \, |y^{-1}|^{d(v)}_v \le |y|_v^{d(v)(m(v) -1)}$.
\end{enumerate}
\end{Lemma}

%%%%%%%%%%%%%%%%%%%%

\begin{proof}
Recall that
\[
|y|_v = \max_{i,j} \{|y^{i,j}(v)|_v\},
\]
where
\[
\rho_v(y) = (y^{i,j}(v))_{i,j} \in \mathrm{M}_{m(v)}(A_v).
\]
Here $|q|_v = |\mathrm{N}_v(q)|^{1/d(v)}_v$ when $q \in A_v$, $\mathrm{N}_v:A_v \to k_v$ is the reduced norm and $\mathrm{dim}_{k_v}(A_v) = d(v)^2$. We noted earlier that $|qq'|_v = |q|_v |q'|_v$ and $|q + q'|_v \le \max\{|q|_v, |q'|_v\}$ for $q, q \in A_v$, so statement \emph{1}.\ of the lemma is clear by the usual matrix multiplication formula.

Let $\lambda_v$ be a prime element of the unique maximal $O_v$-order $U_v$ in $A_v$, so that $\lambda_v U_v = U_v \lambda_v = P_v$ is the maximal two-sided proper ideal of $U_v$. For $0 \ne q \in A_v$ there is an integer $\ell$ such that $q U_v = \lambda_v^\ell U_v$, and $|q|_v = |\mathrm{N}_v(\lambda_v^\ell)|_v^{1 / d(v)} = (\# k(v))^{-\ell / d(v)}$. This interpretation of $|q|_v$ implies that $|a| = \max_{i,j} |a_{i,j}|_v$ is unchanged if we multiply a matrix $a = (a_{i,j}) \in \mathrm{M}_{m(v)}(A_v)$ on the left or right by a permutation matrix, by a matrix which multiplies a single row or column by an element of $U_v^*$, or by an elementary matrix associated with some element of $U_v$. Thus to prove inequality (\emph{2}) of Lemma \ref{lem:detbound1}, we can use these operations to reduce to the case where $y$ is a diagonal matrix with entries $\lambda_v^{z_1},\ldots,\lambda_v^{z_{m(v)}}$ for some integers $z_1, \dots, z_{m(v)}$.

When $y$ has this form,
\[
|y|^{d(v)}_v = \max \{ |\mathrm{N}_v(\lambda_v^{z_i})|_v : 1 \le i \le m(v) \} = (\#k(v))^{-\min_i \{z_i\}},
\]
\[
|\mathrm{det}_v(y)|_v = |\mathrm{N}_v(\lambda)^z|_v = (\# k(v))^{-z} \quad \mathrm{where} \quad z = \sum_{i = 1}^{m(v)} z_i,
\]
and $y^{-1}$ is the diagonal matrix with entries $\lambda_v^{-z_1},\ldots, \lambda_v^{-z_{m(v)}}$.
%Then,
%\[
%|\pi_v^z \lambda_v^{-z_j}|_v = |\mathrm{N}_v(\pi_v^z \lambda_v^{-z_j})^{1/d(v)}|_v =
%\]
%\[
%|\pi_v^{z-z_j}|_v = (\#k(v))^{-(z - z_j)} = (\# k(v))^{z_j - z},
%\]
%where the last $|\ |_v$ is the norm on $k_v$. 
Therefore,
\[
|y^{-1}|^{d(v)}_v = (\# k(v))^{- \min_i \{-z_i\}} = (\# k(v))^{\max_i \{z_i\}}.
\]
The inequality in statement \emph{2}.\ of the lemma is therefore equivalent to
\[
\max_i \{z_i\} - z \le - (m(v)-1) \ \min_i \{z_i\}.
\]
This is the same as
\[
z - \max_i \{z_i\} \ge (m(v) - 1) \min_i \{z_i\},
\]
which is certainly true.
\end{proof}

%%%%%%%%%%%%%%%%%%%%

\begin{Lemma}\label{lem:detbound2}
Suppose $v$ is an infinite place, so that there is an isomorphism $\rho_v : B_v \to \mathrm{M}_{m(v)}(A_v)$ with $A_v = k_v$ if $v$ is complex and either $A_v = k_v$ or $A_v = \mathbb{H}$ if $v$ is real. Define $\mathrm{det}'_v : B_v \to k_v$ by $\mathrm{det}'_v(q) = |\mathrm{det}_v(q)|_v^{1 / d(v)}$ where $\mathrm{det}_v : B_v \to k_v$ is the reduced norm. Then, there are minimal real constants $\delta_1(A_v, m(v))$ and $\delta_2(A_v, m(v))$ such that for all $y, y' \in B_v$:
\begin{enumerate}
\label{eq:upperone2}
\item[1.] $|y y'|_v \le \delta_1(A_v, m(v)) \, |y|_v \, |y'|_v$.
\item[2.] $|y y'|_v = |y|_v \, |y'|_v$ if either $y$ or $y'$ is a scalar matrix or a permutation matrix.
\item[3.] $|\mathrm{det}'_v(y) y^{-1}|_v \le \delta_2(A_v, m(v)) \, |y|_v^{m(v) - 1}$ for all $y \in B_v^*$.
\end{enumerate}
We also have the bounds
\begin{equation}
\label{eq:deltabounder}
1 \le \delta_1(A_v, m(v)) \le m(v)^{[k_v : \mathbb{R}]} 
\end{equation}
\begin{equation}
\label{eq:deltabounder2}
1 \le \delta_2(A_v, m(v)) \le 2^{[k_v : \mathbb{R}] m(v) (m(v) - 1)}.
\end{equation}
Furthermore, if $A_v = k_v$ then
\begin{equation}
\label{eq:deltabounder3}
1 \le \delta_2(A_v, m(v)) \le ((m(v) - 1)!)^{[k_v : \mathbb{R}]}.
\end{equation}
\end{Lemma}

%%%%%%%%%%%%%%%%%%%%

\begin{proof}
As before,
\[
|y|_v = \max_{i,j} \{|y^{i,j}(v)|_v\}
\]
where
\[
\rho_v(y) = (y^{i,j}(v))_{i,j} \in \mathrm{M}_{m(v)}(A_v)
\]
and $|q|_v = |\mathrm{N}_v(q)|^{1 / d(v)}_v$ for $q \in A_v$, and where $\mathrm{N}_v : A_v \to k_v$ is the reduced norm and $\dim_{k_v}(A_v) = d(v)^2$. We noted earlier that $|q q'|_v = |q|_v |q'|_v$ for all $q, q \in A_v$, and by Lemma \ref{lem:triangle-cauchy},
\[
\left| \sum_{i=1}^{m(v)} q_i \right|_v \le m(v)^{[k_v : \mathbb{R}] - 1} \sum_{i=1}^{m(v)} |q_i|_v \le m(v)^{[k_v : \mathbb{R}]} \mathrm{max}_i \{ |q_i|_v \}
\]
for $\{q_i\}_i \subset A_v$. By writing the matrix entries of $y y'$ as sums of products of the entries of $y$ and $y'$ this leads to (\emph{1}) in Lemma \ref{lem:detbound2} and the stated bounds on $\delta_1(A_v, m(v))$. The bound (\emph{2}) in Lemma \ref{lem:detbound2} is  clear.

Now suppose that $y \in B_v^*$. We can find permutation matrices $r$ and $r'$ such that the entry $q$ of $r y r'$ for which $|\ |_v$ is maximal lies in the upper left corner. We then perform Gauss--Jordan elimination on the rows and columns of $r y r'$ to produce matrices $e$ and $e'$ in $\mathrm{M}_{m(v)}(A_v)$ such that $e$ and $e'$ are products of elementary matrices and the nonzero off-diagonal entry of each elementary matrix for $e$ or $e'$ has the form $-\tau / q$ for some entry $\tau$ of $r y r'$, where $|-\tau / q|_v = |\tau|_v / |q|_v \le 1$. The matrix $y_1 = e r y r' e'$ has the same entry $q$ as $y$ in the upper left corner, and all of the other entries in the first row and the first column are $0$. Finally, the other entries of $y_1$ have the form $\alpha - (\tau / q) \beta$ where $\alpha$, $\beta$, and $\tau$ are entries of $r y r'$. Since $|\tau|_v \le |q|_v$ we see from Lemma \ref{lem:triangle-cauchy} that
\[
|\alpha - (\tau / q) \beta|_v \le 2^{[k_v : \mathbb{R}] - 1} (|\alpha|_v + |\beta|_v) \le 2 \cdot 2^{[k_v : \mathbb{R}] - 1} |q|_v = 2^{[k_v : \mathbb{R}]} |y|_v.
\]
Since $q$ is an entry of $y_1$, we deduce that
\[
|y|_v \le |y_1|_v \le 2^{[k_v : \mathbb{R}]} |y|_v \quad \mathrm{and} \quad \det(y_1) = \pm \det(y).
\]
We now continue with $y_1$ and construct matrices $s, s' \in \mathrm{M}_{m(v)}(A_v)$ such that $s$ and $s'$ are products of elementary matrices and permutation matrices, and the off-diagonal entries of the elementary matrices involved in each product have absolute value with respect to $|\ |_v$ bounded above by $1$. The matrix $y' = s y s'$ is diagonal and
\begin{eqnarray}
\label{eq:ytoy1}
|y|_v & \le & |y'|_v \le 2^{[k_v : \mathbb{R}] (m(v) - 1)} |y|_v \\ \det(y') & = & \pm \det(y).
\end{eqnarray}
Then $(y')^{-1} = (s')^{-1} y^{-1} s^{-1}$ and $s' (y')^{-1} s = y^{-1}$, where $s, s', (s')^{-1}$, and $s^{-1}$ are products of elementary matrices and permutation matrices such that the off diagonal entries in each elementary matrix has absolute value with respect to $|\ |_v$ bounded by $1$. This leads by the above reasoning to the bounds
\begin{eqnarray}
\label{eq:ytoy2}
|(y')^{-1}|_v & \le & 2^{[k_v : \mathbb{R}] (m(v) - 1)} |y^{-1}|_v \\
|y^{-1}|_v & \le & 2^{[k_v : \mathbb{R}] (m(v) - 1)} |(y')^{-1}|_v.
\end{eqnarray}

Write $y' = \mathrm{diag}(c_1, \ldots, c_{m(v)})$ for some $c_i \in A_v$. Define $r_i = |c_i|_v = |\mathrm{N}_v(c_i)|_v^{1 / d(v)}$. Then
\[
|y'|_v = \mathrm{max}_i \{r_i\},
\]
\[
\mathrm{det}'(y') = \prod_i r_i = r,
\]
\[
\mathrm{det}'(y') (y')^{-1} = \mathrm{diag}(r c_1^{-1}, \ldots, r c_{m(v)}^{-1}).
\]
We deduce from this that
\begin{eqnarray}
\label{eq:diagb}
|\mathrm{det}'(y') (y')^{-1})|_v & = & \max_i \{|r c_i^{-1}|_v\} \nonumber \\
& = & r \max_i \{|c_i^{-1}|_v\} \nonumber \\
& = & r \max_i \{r_i^{-1}\} \\
& \le & (\max_i \{r_i\})^{m(v) - 1}\\
& = & |y'|_v^{m(v) - 1}.
\end{eqnarray}

Combining this with \eqref{eq:ytoy1} and \eqref{eq:ytoy2} gives
\begin{eqnarray}
\label{eq:finaldiag}
|\mathrm{det}'(y) (y)^{-1})|_v & = & |\mathrm{det}'(y)|_v |y^{-1}|_v \nonumber \\
& = & |\mathrm{det}'(y')|_v |y^{-1}|_v \nonumber \\
& \le & 2^{[k_v : \mathbb{R}] (m(v) - 1)} |\mathrm{det}'(y')|_v |(y')^{-1}|_v \nonumber \\
& \le & 2^{[k_v : \mathbb{R}] (m(v) - 1)} |\mathrm{det}'(y') (y')^{-1}|_v \nonumber \\
& \le & 2^{[k_v : \mathbb{R}] (m(v) - 1)} |y'|_v ^{m(v) - 1} \nonumber \\
& \le & 2^{[k_v : \mathbb{R}] (m(v) - 1)} \left( 2^{[k_v : \mathbb{R}] (m(v) - 1)} |y|_v \right)^{m(v) - 1} \nonumber \\
& = & 2^{[k_v : \mathbb{R}] m(v) (m(v) - 1)} |y|_v^{m(v) - 1}.
\end{eqnarray}
This gives (\emph{2}) in Lemma \ref{eq:upperone2} and the bound \eqref{eq:deltabounder2} on $\delta_2(A_v, m(v))$.

Now, suppose that $A_v = k_v$. We can improve the above bound on $\delta_2(A_v, m(v))$ using the fact that $\det(y) y^{-1}$ is the transpose of the cofactor matrix of $y$. Using the formula for the determinant as a sum over permutations, every entry $\det(y) y^{-1}$ is the sum of $(m(v) -1)!$ terms, each of which is $\pm 1$ times a product of $m(v) - 1$ entries of the matrix $y$. The absolute value with respect to $|\ |_v$ of each entry of $y$ is bounded by $|y|_v$, which implies that
\[
|\det(y) y^{-1}|_v \le ((m(v) - 1)!)^{[k_v : \mathbb{R}] - 1} (m(v) - 1)! |y|_v^{m(v) - 1} =
\]
\[
((m(v) - 1)!)^{[k_v : \mathbb{R}]} |y|_v^{m(v) -1}.
\]
This is the bound in \eqref{eq:deltabounder3}.
\end{proof}

%%%%%%%%%%%%%%%%%%%%

\section{The main result}\label{proof of intro main}

%%%%%%%%%%%%%%%%%%%%

We retain all notation and definitions from $\S\S$\ref{notation}-\ref{sec:heights}. Let $\{\omega_i\}_{i=1}^{n d^2}$ be a $\mathbb{Z}$-basis for $\mathcal{D}$. The discriminant $d_{\mathcal{D}}$ of $\mathcal{D}$ is defined to be 
\[
d_{\mathcal{D}} = \det(M),
\]
where $M$ is the matrix $\left( T(\omega_i \, \omega_j) \right)_{1 \le i,j \le n d^2}$ and $T : B_{\mathbb{R}} \to \mathbb{R}$ is the trace. As a real vector space, $B_{\mathbb{R}} \cong \mathbb{R}^{n d^2}$. The additive Tamagawa measure $\mathrm{Vol}$ on $B$ described in \cite[$\S$X.3]{Cassels--Frohlich} is defined in such a way that
\begin{equation}\label{eq:covol-def}
d_{\mathcal{D}} = \mathrm{Vol}(B_{\mathbb{R}} / \mathcal{D}) = |d_{\mathcal{D}}|^{1/2}.
\end{equation}

Consider a compact convex symmetric subset $X$ of $B_{\mathbb{R}}$. By Minkowski's lattice point theorem, if
\begin{equation}\label{eq:minkowski}
\mathrm{Vol}(X) \geq 2^{\dim_{\mathbb{Q}} B}\ d_{\mathcal{D}},
\end{equation}
then $X$ contains a nonzero element of $\mathcal{D}$. Since $X$ is bounded, there is a constant $m_X$ such that $|\mathrm{Norm}_v(y)|_v$ is bounded by $m_X^{[k_v : \mathbb{R}] / n}$ for every $y \in B_v \cap X$ and $v \in V_\infty$. Then the set
\begin{equation}\label{eq:fund-domain}
F_X = \{ (x, \beta) \in G_S\ :\ x \in X,~\beta \mathcal{D} \subseteq \mathcal{D},~[\mathcal{D} : \beta \mathcal{D}] \leq m_X \}
\end{equation}
is a compact subset of $G_S$.

%%%%%%%%%%%%%%%%%%%%

\begin{Proposition}\label{prop:fund-dom}
With notation as above, suppose that $S$ contains all finite places $v$ of $k$ such that $|\mathrm{Norm}_{k / \mathbb{Q}}(v)|^d \leq m_X$. Then $F_X$ is a fundamental set for the action of $\Gamma_S$ on $G_S$ in the sense that $\Gamma_S F_X = G_S$.
\end{Proposition}

%%%%%%%%%%%%%%%%%%%%

\begin{proof}
Given $(x, \beta) \in G_S$, we must show that there exists $c \in \Gamma_S$ such that $(c x, c \beta) \in F_X$. This happens if and only if
\begin{eqnarray}
c \beta \mathcal{D} & \subseteq & \mathcal{D} \label{eq:cd} \\
\lbrack \mathcal{D} : c \beta \mathcal{D} \rbrack & \leq & m_X, \quad \mathrm{and} \label{eq:Dm} \\
c x & \in & X. \label{eq:cS}
\end{eqnarray}
By definition, \eqref{eq:cd} means that $c \in \mathcal{D} \beta^{-1}$. If $c x \in X$, then
\begin{eqnarray}
\label{eq:bounder!}
[\mathcal{D} : c \beta \mathcal{D}] & = & |\mathrm{Norm}_f(c \beta)|^{-1} \nonumber \\
& = & |\mathrm{Norm}_f(c)|^{-1} \, |\mathrm{Norm}_f (\beta)|^{-1} \nonumber \\
& = & |\mathrm{Norm}_\infty(c)| \, |\mathrm{Norm}_\infty( x)| \\
& = & |\mathrm{Norm}_\infty(cx)| \nonumber \\
& \leq & \prod_{v \in V_\infty} m_X^{[k_v : \mathbb{Q}] / n} = m_X \nonumber
\end{eqnarray}
by \eqref{eq:S-group} and the definitions of $G_S$ and $m_X$. Therefore \eqref{eq:cS} implies \eqref{eq:Dm}. Combining these facts, it suffices to show that $\mathcal{D} \beta^{-1} \cap X x^{-1}$ contains an element of $\Gamma_S$.

Since $X x^{-1}$ is convex and symmetric with volume $\mathrm{Vol}(X) |\mathrm{Norm}_\infty(x)|^{-1}$ and the lattice $\mathcal{D} \beta^{-1}$ in $\mathbb{R}^{n d^2}$ has covolume
\[
\mathrm{Covol}(\mathcal{D} \beta^{-1}) = d_{\mathcal{D}} |\mathrm{Norm}_f(\beta^{-1})|^{-1} = d_{\mathcal{D}} |\mathrm{Norm}_\infty(x)|^{-1},
\]
this implies that
\[
\mathrm{Vol}(X x^{-1}) \geq 2^{\dim_{\mathbb{Q}} B} \mathrm{Covol}(\mathcal{D} \beta^{-1})
\]
if and only if $\mathrm{Vol}(X) \geq 2^{\dim_{\mathbb{Q}} B} d_{\mathcal{D}}$. Since this holds by definition of $X$, it follows that $X x^{-1} \cap \mathcal{D} \beta^{-1}$ contains a nonzero element $c$ of $\mathcal{D} \beta^{-1}$. By construction, $c$ is an element of $B^*$ such that $(c x, c \beta) \in F_X$. We claim that $c \in \Gamma_S$.

Since $c \beta \in \mathcal{D}$, it follows from \eqref{eq:rednorm} that $|\mathrm{Norm}_v((c \beta)_v)|^{-1}$ is a nonnegative integral power of $\mathrm{Norm}_{k / \mathbb{Q}}(v)^d$ for each $v \in V_f$. We know that
\[
\prod_{v \not \in V_\infty} |\mathrm{Norm}_v((c \beta)_v)|^{-1} = |\mathrm{Norm}_f(c \beta)|^{-1} = |\mathrm{Norm}_\infty(c x)| \leq m_X
\]
by \eqref{eq:bounder!}. Hence if $|\mathrm{Norm}_v((c \beta)_v)| \ne 1$ for some finite place $v$, then
\[
|\mathrm{Norm}_{k / \mathbb{Q}}(v)|^d \leq m_X,
\]
which implies that $v \in S_f$. It follows that
\[
|\mathrm{Norm}_v((c \beta)_v)| = 1
\]
for all $v \in V_f \smallsetminus S_f$. Thus $(c\beta)_v \mathcal{D}_v = \mathcal{D}_v$ for these $v$, since $c \beta \in \mathcal{D}$. However, $\beta_v = 1$ if $v \notin S$, so
\[
c \mathcal{D}_v = c_v \mathcal{D}_v = (c\beta)_v \mathcal{D}_v = \mathcal{D}_v
\]
for all $v \in V_f \smallsetminus S_f$. This implies that $c \in \mathcal{D}^*_v$ for all $v \notin S_f$, so $c \in \Gamma_S$. This proves the proposition.
\end{proof}

%%%%%%%%%%%%%%%%%%%%

We now describe how $F_X$ determines generators for $\Gamma_S$. A subset $P$ of $G_S$ will be called a set of \emph{topological generators} for $G_S$ if for any open subset $O$ of $G_S$, the group generated by $O$ and $P$ is all of $G_S$. The following lemma should be compared with \cite[Lemma 6.3]{Lenstra}.

%%%%%%%%%%%%%%%%%%%%

\begin{Lemma}\label{lem:top-gens}
Let $P$ be a set of topological generators for $G_S$ that contains the identity, and let $F_X$ be as in Proposition \ref{prop:fund-dom}. Then $\Gamma_S$ is generated by its intersection with $F_X P F_X^{-1}$.
\end{Lemma}

%%%%%%%%%%%%%%%%%%%%

\begin{proof}
We have an equality of sets
\[
F_X (P \cup P^{-1}) F_X^{-1} = (F_X P F_X^{-1}) \cup (F_X P F_X^{-1})^{-1}.
\]
Therefore we can replace $P$ by $P \cup P^{-1}$ for the remainder of the proof and assume that $P$ is symmetric, i.e., that $P = P^{-1}$. We emphasize that this does not mean we must assume $P$ is symmetric in the statement of the lemma.

Consider the subset
\[
O = (G_S \smallsetminus \Gamma_S) \cup (\Gamma_S \cap F_X F_X^{-1})
\]
of $G_S$. This is an open neighborhood of $F_X F_X^{-1}$ in $G_S$ because $\Gamma_S$ is a discrete subgroup of $G_S$. Since $F_X$ is a fundamental set for the action of $\Gamma_S$ on $G_S$, we can find a subset $F \subset F_X$ such that $\Gamma_S \times F \to G_S$ is a bijection. We claim that there is a small open neighborhood $U$ of the identity in $G_S$ such that $F U F^{-1} \subset O$.

It will be enough to find a $U$ such that $\Gamma_S \cap (F U F^{-1}) \subset F_X F_X^{-1}$. Let $T$ be the set of $\gamma \in \Gamma_S$ such that $\gamma F \cap F U \ne \emptyset$. We want to show that if $\gamma \in T$, then $\gamma F_X \cap F_X \ne \emptyset$. Since $F$ is a bounded fundamental domain for the action of $\Gamma_S$ on $G_S$ and $\Gamma_S$ is discrete in $G_S$, the set $T$ is finite when $U$ is bounded. We can then shrink $U$ further and assume that if $\gamma \in T$, then $\gamma F' \cap F' \ne \emptyset$ for each open neighborhood $F'$ of the closure of $F$ in $G_S$. If $\gamma F_X \cap F_X = \emptyset$ for some $\gamma \in T$, then since $F_X$ is compact there will be an open neighborhood $F'$ of $F_X$ for which $\gamma F' \cap F' = \emptyset$. This is a contradiction, since the closure of $F$ is contained in $F_X$, which proves the claim.

Let $P^\prime = P \cup U$, so $\langle P^\prime \rangle = G_S$, and let $\Delta < \Gamma_S$ be the subgroup generated by $\Gamma_S \cap F P^\prime F^{-1}$. We claim that $\Delta = \Gamma_S$. Indeed, if $x p \in F P^\prime$, there exist $y \in F$ and $\gamma \in \Gamma_S$ such that $xp = \gamma y$. Then
\[
\gamma = x p y^{-1} \in F P^\prime F^{-1},
\]
so $\gamma \in \Delta$. This implies that $F P^\prime \subseteq \Delta F$, so $\Delta F P^\prime \subseteq \Delta F$. Therefore, $\Delta F$ is right $P^\prime$-invariant, but $\langle P^\prime \rangle = G_S$, so $\Delta F = G_S$. Since $\Gamma_S \times F \to G_S$ is a bijection, it follows that $\Delta = \Gamma_S$.

This proves that $\Gamma_S$ is generated by
\[
\Gamma_S \cap F P^\prime F^{-1} \subseteq (\Gamma_S \cap F P F^{-1}) \cup (\Gamma_S \cap F U F^{-1}).
\]
However, $F U F^{-1} \subset O$, and $\Gamma_S \cap O \subset F_X F_X^{-1}$ by definition, so
\begin{equation}\label{eq:gen-set}
\Gamma_S \cap F P^\prime F^{-1} \subseteq (\Gamma_S \cap F_X P F_X^{-1}) \cup (\Gamma_S \cap F_X F_X^{-1}).
\end{equation}
Since $P$ contains the identity, the right side of \eqref{eq:gen-set} equals $\Gamma_S \cap F_X P F_X^{-1}$. This proves the lemma.
\end{proof}

%%%%%%%%%%%%%%%%%%%%

We now define several constants that we need to state our main result.
\begin{enumerate}

\item For $X$, $F_X$, and $\ell$ as above, let $T_1$ be the supremum of $1$ and
\[
\left\{ |x_v|_v^{d(v) / [k_v : \mathbb{R}])} \right\}
\]
over all
\[
x = \prod_{v \in V_\infty} x_v \in B_{\mathbb{R}}
\]
for which $(x,\beta) \in F_X$ for some $\beta$.

\item Let $P$ be a finite set of topological generators for $G_S$ which contains the identity element (see \S \ref{s:topgens} for an example of such a set). We assume that every element of $P$ has the form $(z, \zeta)$ with $z = \prod_{v \in S_\infty} z_v \in B_{\mathbb{R}}^*$ and $\zeta = \prod_{v \in S_f} \zeta_v \in B_{S_f}^*$, where $z_v$ is a real scalar and each $\zeta_v$ lies in the local maximal order $\mathrm{M}_{m(v)}(U_v)$ of $B_v = \mathrm{M}_{m(v)}(A_v)$ (cf.\ $\S$\ref{s:topgens}). Let $T_2$ be the supremum of $1$ and 
\[
\left \{ |z_v|_v^{d(v)/ [k_v : \mathbb{R}]} \right \}
\]
over all $z = \prod_v z_v \in P$ and all $v \in V_\infty$.

\item Let $T_3$ be
\[
\prod_{v \in S_f} \max \left\{ 1, |\alpha_v|^{d(v)}_v \right\}
\]
where as $(x,\alpha)$ ranges over $F_X$ and $\alpha = \prod_{v \in S_f} \alpha_v \in B_{S_f}^*$. Note that for such $\alpha$ and $\alpha_v$ we have that $\alpha_v \mathcal{D}_v \subseteq \mathcal{D}_v$. Such $\alpha_v $ are contained in $\mathcal{D}_v$, so this constant is finite. Similarly, define $T_3'$ to be the maximum of
\[
\prod_{v \in S_f} \max \left\{ 1, |\alpha_v|_v^{d(v)(m(v) - 1)} \right\},
\]
where $\alpha$ and the $\alpha_v$ range as above.

\item Let $T_4$ be the smallest number such that
\[
\prod_{v \in S_f} \max \left\{ 1, |g_v|^{d(v)}_v \right\} \leq T_4
\]
for all $g = \prod_v g_v \in P$.

\item Let $T_5$ be the supremum of $1$ and
\[
\left\{ |\mathrm{det}_v(a_v)|_v^{1 / [k_v : \mathbb{R}]}\ :\ a \in F_X\ \textrm{and}\ v \in V_\infty \right\}.
\]
where we write $a \in F_X$ as $a = (a_v)_v$ with $a_v \in B_v$. (Recall that $\mathrm{det}_v:B_v \to k_v$ is the reduced norm.)

\item Let $T_6$ be the maximum over all subsets $W$ of $V_\infty$ of
\[
T_1^{a(W)} \, T_2^{b(W)} \, T_5^{b(V_\infty \smallsetminus W)}
\]
where
\begin{equation}
\label{eq:abdef}
a(W) = \sum_{v \in W} [k_v : \mathbb{R}] m(v) \quad \mathrm{and} \quad b(W) = \sum_{v \in W} [k_v : \mathbb{R}].
\end{equation}

\end{enumerate}
Now we are ready to state and prove our main result.

%%%%%%%%%%%%%%%%%%%%

\begin{Theorem}\label{bounded height}
Let $k$ be an algebraic number field of degree $n$ over $\mathbb{Q}$ and $B$ a central simple $k$-division algebra of degree $d$. Let $S$ be a finite set of places of $k$ containing all the archimedean places $V_\infty$ and let $\mathcal{D} \subset B$ be an $O_k$-order. We suppose that the $k_v$ isomorphisms $\rho_v : B_v \to \mathrm{M}_{m(v)}(A_v)$ are chosen such that $\rho_v(\mathcal{D}_v) \subseteq \mathrm{M}_{m(v)}(U_v)$ for $v \not \in S$, where $U_v$ is the unique maximal $O_v$-order in the $k_v$-division algebra $A_v$. Suppose $s$ is the number of (real) places $v$ at which $A_v$ is isomorphic to $\mathbb{H}$. Let $G_S$ be the topological group defined in \eqref{eq:S-group}, $P$ be a topological generating set for $G_S$ satisfying the above conditions, and let $\Gamma_S$ be the group of $S$-units associated with $\mathcal{D}$.

Suppose that $X$ is a convex symmetric subset of $B_{\mathbb{R}}$ such that \eqref{eq:minkowski} holds, and let $m_X$ be the smallest real number such that $|\mathrm{Norm}_v(y)|_v$ is bounded by $m_X^{[k_v : \mathbb{R}] / n}$ for every $y \in B_v \cap X$ and $v \in V_\infty$. Suppose that $S$ contains every finite place $v$ of $k$ such that $\mathrm{Norm}(v) \le m_X^{1/d}$. Finally, let $T_1, \dots, T_6$ be the constants defined immediately above.

Then the set $\Gamma_S \cap F_X P F_X^{-1}$ from Lemma \ref{lem:top-gens} is contained in
\begin{equation}
\label{eq:target}
\mathfrak G_{S, X} = \mathrm{BH}_{S} \left( \left( (d - 1)! d\right)^n \left( \frac{2^{(d / 2)(d - 2)} d}{4 (d - 1)!} \right)^s T_6 T_3 T'_3 T_4 \right).
\end{equation}
Consequently, $\mathfrak G_{S, X}$ is a finite generating set for $\Gamma_S$.
\end{Theorem}

%%%%%%%%%%%%%%%%%%%%

\begin{proof}
Suppose that $\gamma \in \Gamma_S \cap F_X P F_X^{-1}$. Then, there exist elements $(z, \zeta) \in P$ and
\[
(x, \alpha), (y, \beta) \in F_X, \quad x, y \in X, \quad \alpha = \prod_{v \in S_f} \alpha_v,\ \beta = \prod_{v \in S_f} \beta_v
\]
so that $(\gamma, \gamma) = (x, \alpha) (z, \zeta) (y^{-1}, \beta^{-1})$. That is, $x_v z_v y_v^{-1} = \gamma$ for each $v \in V_\infty$ and $\alpha_v \zeta_v \beta_v^{-1} = \gamma$ for each $v \in S_f$.

Let $W(\gamma) = W_\infty(\gamma) \cup W_f(\gamma)$ be the set of places $v$ of $k$ at which $|\gamma|_v > 1$. By assumption, if $v \not \in S$, then $v$ is finite and $\rho_v(\mathcal{D}_v) \subseteq \mathrm{M}_{m(v)}(U_v)$. Thus $\gamma \in \Gamma_S$ implies that $|\gamma|_v \le 1$ if $v \not \in S$. Thus $W(\gamma) \subseteq S$, $W_\infty(\gamma) \subset V_\infty$ and $W_f(\gamma) \subseteq S_f$.

By definition of $H(\gamma)$ we have 
\[
H(\gamma) = \prod_{v \in W_\infty( \gamma)} |x_v z_v y_v^{-1}|^{d(v)}_v \times \prod_{v \in W_f( \gamma)} \mathrm |\alpha_v \zeta_v \beta_v^{-1}|^{d(v)}_v,
\]
Recall that $\mathrm{det}_v:B_v \to k_v$ is the reduced norm and that $d(v)^2$ is the dimension of $A_v$ over $k_v$. If $v$ is archimedean, we defined $\mathrm{det}'_v : B_v \to k_v$ by $\mathrm{det}'_v(q) = |\mathrm{det}_v(q)|_v^{1/d(v)}$ for $q \in B_v$. 
%For $v$ nonarchimedean, $\mathrm{det}'_v : B_v \to k_v$ was defined so that $\mathrm{det}'_v(q)^{d(v)}$ generates the same $O_v$ ideal as $\mathrm{det}_v(q)$ for any $q \in B_v$.

We have $|c \alpha|_v = |c|_v |\alpha|_v$ for $c \in k_v$ and $\alpha \in B_v$, where $|c|_v$ here denotes the absolute value of $c \in k_v$ with respect to $|\ |_v : k_v \to \mathbb{R}$. Therefore we can rewrite the above expression as
\begin{eqnarray}
\label{eq:bigeq}
H(\gamma) & = & \nonumber\\
& & \prod_{v \in W_\infty(\gamma)} |\mathrm{det}'_v(y_v) x_v z_v y_v^{-1}|^{d(v)}_v \label{eq:bound1} \\
& \times & \prod_{v \in W_f( \gamma)} \left ( |\mathrm{det}_v(\beta_v)|_v \, | \alpha_v \zeta_v \beta_v^{-1}|^{d(v)}_v \right )
\label{eq:bound2} \\
& \times & \prod_{v \in W_\infty(\ \gamma)} |\mathrm{det}_v(y_v)|_v^{-1} \label{eq:bound3} \\
& \times & \prod_{v \in W_f(\gamma)} |\mathrm{det}_v(\beta_v)|_v^{-1} \label{eq:bound4},
\end{eqnarray}
where the last two products are computed using the absolute values on the completions $k_v$. We now proceed to bound each of these terms.

%%%%%%%%%%%%%%%%%%%%

\begin{Sublemma}\label{first term bound}
One has that
\begin{equation}\label{eq:sublemma1eq}
\prod_{v \in W_\infty(\gamma)} |\mathrm{det}'_v(y_v) (x_v z_v y_v^{-1})|^{d(v)}_v \leq \mu_1 \ T_1^{a(W_\infty(\gamma))}\ T_2^{b(W_\infty(\gamma))}
%&\le& \mu_1 \left(T_1^{d} T_2 \right)^{ [k:\mathbb{Q}]},
\end{equation}
where $a(W_\infty(\gamma))$ and $b(W_\infty(\gamma))$ are as in (\ref{eq:abdef}) and
\begin{equation}
\label{eq:delbound}
\mu_1 = \prod_{v \in W_\infty(\gamma)} \delta_1(A_v, m(v))^{d(v)} \delta_2(A_v, m(v))^{d(v)}.
\end{equation}
Furthermore,
\begin{equation}
\label{eq:delbound2}
\mu_1 \le ( (d - 1)!d)^{[k : \mathbb{Q}]} \left( \frac{ 2^{(d / 2)(d - 2)} d}{4 (d - 1)!} \right)^s.
\end{equation}
if $A_v = \mathbb{H}$ at exactly $s$ real places of $k$.
\end{Sublemma}

%%%%%%%%%%%%%%%%%%%%

\begin{proof}
By assumption each $z_v$ is a real scalar. Therefore, Lemma \ref{lem:detbound2} gives
\begin{eqnarray}
\label{eq:easy}
& & |\mathrm{det}'_v(y_v) (x_v z_v y_v^{-1})|^{d(v)}_v \nonumber \\
& \le & \delta_1(A_v, m(v))^{d(v)} \, |x_v|^{d(v)}_v \, |z_v|^{d(v)}_v \, |\mathrm{det}'_v(y_v) y_v^{-1}|_v^{d(v)} \\
& \le & \delta_1(A_v, m(v))^{d(v)} \, |x_v|^{d(v)}_v \, |z_v|^{d(v)}_v \, \delta_2(A_v, m(v))^{d(v)} |y_v|_v^{d(v)(m(v) - 1)} \nonumber \\
& \le & \delta_1(A_v, m(v))^{d(v)} \, \delta_2(A_v, m(v))^{d(v)} \, T_1^{[k_v : \mathbb{R}]} \, T_2^{[k_v : \mathbb{R}]} \, T_1^{[k_v : \mathbb{R}] (m(v) - 1)} \nonumber
\end{eqnarray}
Taking the product over all $v \in W_\infty(\gamma)$, we get \eqref{eq:sublemma1eq} and \eqref{eq:delbound}
because $T_1, T_2 \ge 1$. To prove the bound in \eqref{eq:delbound2} we first note that $(m(v), d(v)) = (d, 1)$ if $v$ is archimedean and $A_v = k_v$, while $(m(v), d(v)) = (d / 2, 2)$ if $k_v = \mathbb{R}$ and $A_v = \mathbb{H}$. By Lemma \ref{lem:detbound2},
\[
\delta_1(A_v, m(v))^{d(v)} \delta_2(A_v, m(v))^{d(v)} \le d^{[k_v : \mathbb{R}]} (d - 1)!^{[k_v : \mathbb{R}]}
\]
if $A_v = k_v$ and
\[
\delta_1(A_v, m(v))^{d(v)} \delta_2(A_v, m(v))^{d(v)} \le (d / 2)^{2 [k_v : \mathbb{R}]} 2^{[k_v : \mathbb{R}] (d / 2) (d - 2)}
\]
when $A_v = \mathbb{H}$. Since $\delta_1(A_v, m(v)) \ge 1$ and $\delta_2(A_v, m(v)) \ge 1$ for all archimedean $v$, and $[k_v : \mathbb{R}] = 1$ if $A_v = \mathbb{H}$ we see that
\begin{eqnarray}
\mu_1 & = & \prod_{v \in W_\infty(\gamma)} \delta_1(A_v, m(v))^{d(v)} \delta_2(A_v, m(v))^{d(v)} \nonumber \\
& \le & \prod_{v \in V_\infty} \delta_1(A_v, m(v))^{d(v)} \delta_2(A_v, m(v))^{d(v)} \nonumber \\
& \le & \left( \prod_{v \in V_\infty} ((d - 1)!d)^{[k_v : \mathbb{R}]} \right) \left( \frac{(d / 2)^2 2^{(d / 2) (d - 2)}}{(d - 1)!d} \right)^s.
\end{eqnarray}
This gives \eqref{eq:delbound2} since $\sum_{v \in V_\infty} [k_v : \mathbb{R}] = n$.
\end{proof}

%%%%%%%%%%%%%%%%%%%%

\begin{Sublemma}\label{second term bound}
\begin{equation}
\label{eq:oyboy}
\prod_{v \in W_f(\gamma)} \left ( |\mathrm{det}_v(\beta_v) |_v \, |(\alpha_v \zeta_v \beta_v^{-1})|^{d(v)}_v \right ) \leq T_3  T_4 T'_3 \le T_3^d T_4.
\end{equation}
\end{Sublemma}

%%%%%%%%%%%%%%%%%%%%

\begin{proof}
By Lemma \ref{lem:detbound1}, for every $v \in W_f(\gamma)$ we have that
\begin{eqnarray}
\label{eq:easy2}
|\mathrm{det}_v(\beta_v) |_v \, | (\alpha_v \zeta_v \beta_v^{-1})|^{d(v)}_v & \le & |\alpha_v|^{d(v)} \, |\zeta_v|^{d(v)} \, |\mathrm{det}_v(\beta_v)|_v \, |\beta_v^{-1}|^{d(v)}_v \nonumber \\
& \le & |\alpha_v|^{d(v)}_v \, |\zeta_v|^{d(v)}_v \, | \beta_v|_v^{d(v)(m(v) - 1)}.
\end{eqnarray}
We take the product of these bounds over $v \in W_f(\gamma)$ to deduce \eqref{eq:oyboy}.
\end{proof}

%%%%%%%%%%%%%%%%%%%%

\begin{Sublemma}\label{third term bound}
\[
\prod_{v \in W_\infty(\gamma)} |\mathrm{det}_v(y_v)|_v^{-1} \leq |\mathrm{Norm}_\infty(y)|^{-1 / d} \, T_5^{b(V_\infty \smallsetminus W_\infty(\gamma))}.
\]
\end{Sublemma}

%%%%%%%%%%%%%%%%%%%%

\begin{proof}
We have 
\[
\prod_{v \in W_\infty(\gamma)} |\mathrm{det}_v(y_v)|_v^{-1} =
\]
\[
\frac{1}{|\mathrm{Norm}_\infty(y)|^{1 / d}} \prod_{v \in V_\infty  \smallsetminus  W_\infty(\gamma)} |\mathrm{det}_v(y_v)|_v.
\]
since $\mathrm{Norm}_\infty$ is associated with the $d^{th}$ power of the reduced norm. %Since
%\[
%|\mathrm{Norm}_\infty(y)| = |\mathrm{Norm}_f(\beta)|^{-1} = [\mathcal D : \beta \mathcal D],
%\]
%we see that
%\[
%\frac{1}{|\mathrm{Norm}_\infty(y)|^{1/d}} \leq 1.
%\]
Since $y \in F_X$, by the definition of $T_5$ we have $|\det(y_v)|_v \le T_5^{[k_v : \mathbb{R}]}$ for all $v \in V_\infty$, so the lemma is clear.
% Therefore
%\[
%\prod_{v \in W_\infty( \gamma)} |\mathrm{det}_v(y_v)|_v^{-1} \le \prod_{v \in V_{\infty} \smallsetminus W_\infty(\gamma)} |\mathrm{det}_v(y_v)|_v
%\]
%\[
%\leq \prod_{v \in V_\infty} T_5^{[k_v : \mathbb R]} = T_5^{[k : \mathbb Q]}.
%\]
%This bounds \eqref{eq:bound3}.
\end{proof}

%%%%%%%%%%%%%%%%%%%%

\begin{Sublemma}\label{fourth term bound}
\[
\prod_{v \in W_f(\gamma)} |\mathrm{det}_v(\beta_v)|_v^{-1} \leq |\mathrm{Norm}_f(\beta)|^{-1 / d} = |\mathrm{Norm}_\infty(y)|^{1 / d}
\]
\end{Sublemma}

%%%%%%%%%%%%%%%%%%%%

\begin{proof} Recall that $(y,\beta) \in F_X$, so that $\beta \mathcal{D} \subset \mathcal{D}$.  Since $\mathcal{D}$
is a $O_k$-lattice in $B$, this implies that all the components of $\beta = \prod_{v \in S_f} \beta_v \in B_f^*$
are integral over $O_k$.  We view $\beta$ an idele in $J(B)$ with component $1$ outside of $S_f$.  Since $\mathrm{det}_v:B_v \to k_v$
is the reduced norm, we conclude that $|\mathrm{det}_v(\beta_v)|_v \le 1$ for all places $v$ of $k$. 
We now have
\[
\prod_{v \in W_f( \gamma)} |\mathrm{det}_v(\beta_v)|_v^{-1} =
\]
\[
\prod_{v \in V_f} |\mathrm{det}_v(\beta_v)|_v^{-1} \, \prod_{v \in V_f \smallsetminus  W_f(\gamma)} |\mathrm{det}_v(\beta_v)|_v
\]
\[ 
\le \prod_{v \in V_f} |\mathrm{det}_v(\beta_v)|_v^{-1} = |\mathrm{Norm}_f(\beta)|^{-1 / d}.
\]
The second equality in the statement of the sublemma follows from the definition of $(y,\beta) \in G_S$.
\end{proof}

%%%%%%%%%%%%%%%%%%%%

Substituting these bounds in \eqref{eq:bound1} - \eqref{eq:bound4} shows that $\gamma$ lies in \eqref{eq:target} and completes the proof of Theorem \ref{bounded height}.
\end{proof}

%%%%%%%%%%%%%%%%%%%%

\begin{Remark}
Notice that we used the fact that $B$ is a division algebra to ensure that Minkowski's theorem returns an invertible element of $B$. It would be interesting to prove an analogue of our result for $S$-unit groups of any central simple algebra over a number field. Given that $\mathrm{SL}_n(\mathbb{Z})$ is generated by elementary matrices, which certainly have small height, we expect such a result to hold.
\end{Remark}

%%%%%%%%%%%%%%%%%%%%

\section{Explicit bounds}

%%%%%%%%%%%%%%%%%%%%

In this section we make some particular choices in order to give more explicit calculations of the bounds in the previous section.

As before, $B$ is a division algebra of dimension $d^2$ with center a number field $k$ of degree $n = [k:\mathbb{Q}]$ over $\mathbb{Q}$. Let $V = V_\infty \cup V_f$ be the set of places of $k$, $r_1$ the number of real places of $k$, and $r_2$ the number of complex places. For $v \in V$, we fix an isomorphism of $B_v$ with $\mathrm{M}_{m(v)}(A_v)$ where $A_v$ is a division algebra of dimension $d(v)^2$ over its center $k_v$.

%%%%%%%%%%%%%%%%%%%%

\subsection{A maximal order $\mathcal{D}$ and an archimedean set $X$}

%%%%%%%%%%%%%%%%%%%%

Suppose that $\mathcal{D}$ is the maximal order of $B$ that is isomorphic to $\mathrm{M}_{m(v)}(U_v)$ for all finite $v$, where $U_v$ is the unique maximal $O_v$ order in $A_v$.

Suppose first that $v$ is archimedean. The normalized Haar measure on $k_v$ is the Euclidean measure if $k_v = \mathbb{R}$ and is twice the Euclidean measure if $k_v = \mathbb{C}$. The normalized Haar measure on $\mathbb{H} = \mathbb{R} + \mathbb{R} I + \mathbb{R} J + \mathbb{R} I J$ is $4$ times the one associated to the usual Euclidean measure on $\mathbb{R}^4$ under the basis $\{1, I, J, IJ\}$. The norm on $\mathbb{H}$ is $|\alpha| = |N(\alpha)|^{1 / 2}$ where $N : \mathbb{H} \to \mathbb{R}$ is the reduced norm. Thus the volume with respect to the normalized Haar measure of a ball of radius $c^{1 / 2}$ inside $\mathbb{H}$ is
\[
4 \frac{\pi^2}{\Gamma(3)} c^2 = 2 \pi^2 c^2.
\]
The normalized volume of a ball of radius $c$ in $\mathbb{C}$ (resp.\ $\mathbb{R}$) is $2 \pi c^2$ (resp.\ $2 c$). The normalized Haar measure on $B_v = \mathrm{M}_{m(v)}(A_v)$ is then the product measure associated with matrix entries. Here $(m(v), d(v)) = (d, 1)$ if $A_v = k_v$ and $(m(v), d(v))= (d / 2, 2)$ if $A_v = \mathbb{H}$, and there are $m(v)^2$ matrix entries associated to each element of $B_v$.

Let $S_{ram, \infty}(B)$ be the set of infinite places of $k$ at which $B$ ramifies. Recall that $s = \#S_{ram,\infty}(B)$. We have $\dim_{k_v}(A(v)) = d(v)^2$, so $d(v) = 1$ or $2$ if $v$ is infinite. Let $c > 1$ be a real parameter, and let $X(c)$ be the set of
\[
x = \prod_{v \in V_\infty} x_v \in B_{\mathbb{R}} = \prod_{v \in V_\infty} B_v
\]
such that $|x|_v^{d(v) / [k_v : \mathbb{R}]} \le c$ for all $v \in V_\infty$. Then
\[
\mathrm{Vol}(X(c)) = (2 c)^{d^2 (r_1 - s)} \left( 2 \pi^2 c^2 \right)^{(d / 2)^2 s} (2 \pi c^2)^{d^2 r_2} = z c^{d^2(n - s / 2)}
\]
where 
\begin{equation}
\label{eq:zdef}
z = 2^{d^2 (r_1 - s)} \left( 2 \pi^2 \right)^{(d / 2)^2 s} (2 \pi )^{d^2 r_2}.
\end{equation}
Now choose $c$ such that
\begin{equation}\label{eq:minkowski2}
z c^{d^2 (n - s / 2)} = \mathrm{Vol}(X(c)) = 2^{\dim_{\mathbb{Q}} B}\ d_{\mathcal{D}} = 2^{d^2 n}\ d_{\mathcal{D}}.
\end{equation}
In other words,
\begin{equation}
\label{eq:cdef}
c = \left( \frac{2^{d^2 n}}{z} d_{\mathcal{D}} \right)^{\frac{1}{d^2 (n - s / 2)}} = \left( \frac{2}{\pi} \right)^{\frac{r_2}{n - s / 2}} \left( \frac{2 \sqrt{2}}{\pi} \right)^{\frac{s}{2 n - s}} d_{\mathcal{D}}^{\frac{1}{d^2 (n - s / 2)}}.
\end{equation}

%%%%%%%%%%%%%%%%%%%%

\subsection{The constant $m_X$}

%%%%%%%%%%%%%%%%%%%%

Setting $X = X(c)$, we need to find an $m_X$ such that
\begin{equation}
\label{eq:normreq}
|\mathrm{Norm}_v(y_v)|_v = |\mathrm{det}_v(y_v)|^{d} \le m_X^{[k_v : \mathbb{R}] / n}
\end{equation}
for all $v \in V_\infty$, where
\[
y = \prod_{v \in V_\infty} y_v \in X(c).
\]
and $\mathrm{det}_v:B_v \to k_v$ is the reduced norm.  

%%%%%%%%%%%%%%%%%%%%
 
\subsubsection{Real $v$ with $A_v = k_v$}

%%%%%%%%%%%%%%%%%%%%
 
In this case $\mathrm{det}_v(y_v)$ is the determinant of a real $d \times d$ matrix each of whose entries is bounded by $c$ in absolute value. The Euclidean length of each column of $y_v$ is thus bounded by $d c$, so we have
\begin{equation}
\label{eq:realbound}
|\mathrm{Norm}(y_v)|_v = |\mathrm{det}_v(y_v)|_v^d  \le ((c d)^d)^d = (c d)^{d^2}
\end{equation}
whenever $k_v = A_v = \mathbb{R}$.

%%%%%%%%%%%%%%%%%%%%

\subsubsection{Real $v$ with $A_v = \mathbb{H}$}

%%%%%%%%%%%%%%%%%%%%

We have a representation
\begin{equation}
\label{eq:Hrep} 
\mathbb{H} = \mathbb{R} + \mathbb{R} I + \mathbb{R} J + \mathbb{R} I J \to \mathrm{M}_2(\mathbb{C})
\end{equation}
determined by
\[
I \to \begin{pmatrix} \sqrt{-1} & 0 \\ 0 & -\sqrt{-1} \end{pmatrix}, \quad J \to \begin{pmatrix} 0 & 1 \\ -1 & 0 \end{pmatrix}, \quad I J \to \begin{pmatrix} 0 & \sqrt{-1} \\ \sqrt{-1} & 0 \end{pmatrix}.
\]
This sends $y_v \in \mathrm{M}_{d / 2}(\mathbb{H})$ to a matrix $y'_v \in \mathrm{M}_{d}(\mathbb{C})$ which consists of $2 \times 2$ blocks of the form
\[
\begin{pmatrix} \lambda & \mu \\
-\overline{\mu} & \overline{\lambda} \end{pmatrix}
\]
for $\lambda, \mu \in \mathbb{C}$.

Since $y_v \in X(c)$, we must have $|y_v|_v^{d(v)} = |y_v|_v^2 \le c$, where $|y_v|_v^2$ is the supremum of $|\lambda|^2 + |\mu|^2$ over the above $2 \times 2$ blocks. The columns of $y'_v$ are vectors $y'_v(1),\ldots, y'_v(d)$ in $\mathbb{C}^{d}$ with the property that with respect to the usual Hermitian inner product $\langle \ , \ \rangle$ on $\mathbb{C}^d$ we have
\[
\langle y'_v(i), y'_v(i) \rangle \le d c / 2 \quad \mathrm{and} \quad \langle y'_v(2 j - 1),y'_v(2 j) \rangle = 0
\]
for $1 \le i \le d$ and $1 \le j \le d / 2$.

Let $T$ be the subset of $\mathbb{C}^d$ consisting of all linear combinations of the form $\sum_{i = 1}^d \tau_i y'_v(i)$, where $|\tau_i|_{\mathbb{C}} \le 1$ and $|\ |_{\mathbb{C}}$ is the usual Euclidean inner product on $\mathbb{C}$. Recall that the normalized Haar measure on $\mathbb{C}$ is $2$ times the standard Euclidean Haar measure. Give $\mathbb{C}^d$ the product measure. Then $T$ is the image of the unit polydisc in $\mathbb{C}^d$ under left multiplication by the matrix $y'_v$. Therefore
\begin{equation}
\label{eq:volt1}
\mathrm{Vol}(T) = |\det(y')|^2_{\mathbb{C}} (2 \pi)^d.
\end{equation}
On the other hand, $T$ is contained in the product of real two-dimensional metric disks $T_i$ defined by 
\[
T_i = \{ \tau_i \ y'_v(i) : |\tau_i| \le 1 \}
\]
as $i$ ranges over $1 \le i \le d$.  The square Euclidean length $\langle y'_v(i), y'_v(i) \rangle$ of each $y'_v(i)$ is bounded by $d c / 2.$ We conclude from this that
\begin{equation}
\label{eq:volt2}
\mathrm{Vol}(T) \le (2 \pi dc/2)^d.
\end{equation}
The reduced norm $\mathrm{det}_v(y_v)$ is equal to $\det(y'_v)$, so \eqref{eq:volt1} and \eqref{eq:volt2} give
\begin{equation}
\label{eq:volt3}
|\mathrm{Norm}_v(y_v)|_v = |\mathrm{det}_v(y_v)|_v^d = |\det(y')|^d_{\mathbb{C}} \le (d c / 2)^{d^2 / 2}
\end{equation}
when $k_v = \mathbb{R}$ and $A_v = \mathbb{H}$.

%%%%%%%%%%%%%%%%%%%%

\subsubsection{Complex $v$}

%%%%%%%%%%%%%%%%%%%%

Finally suppose $k_v = A_v = \mathbb{C}$. We can define $y'_v = y_v$ and use the above arguments to bound $|\mathrm{Norm}_v(y_v)|_v = |\det(y'_v)|^{2 d}_{\mathbb{C}}$. Here the columns of $y'_v$ have complex square length bounded by $d c^2$ since $|y_v|_v^{1 / 2} \le c$ in this case, where $|y_v|_v$ is the supremum of the normalized absolute value with respect to $|\ |_v$ of the matrix coefficients of $y_v = y'_v$, and $|\ |_v$ is $|\ |_{\mathbb{C}}^2$. Using the same set $T$ defined above leads to
\[
\mathrm{Vol}(T) = |\det(y')|^2_{\mathbb{C}} (2 \pi)^d \le (2 \pi dc^2)^d.
\]
We see from this that
\begin{equation}
\label{eq:volt4}
|\mathrm{Norm}_v(y_v)|_v = |\det(y')|^{2 d}_{\mathbb{C}} \le (d c^2)^{d^2}
\end{equation}
when $k_v = \mathbb{C} = A_v$.

%%%%%%%%%%%%%%%%%%%%

\subsubsection{A choice for $m_X$}

%%%%%%%%%%%%%%%%%%%%

We can now put together \eqref{eq:volt2}, \eqref{eq:volt3}, and \eqref{eq:volt4} to find a constant $m_X$ that satisfies \eqref{eq:normreq}. One can take
\begin{eqnarray}
\label{eq:mxdef}
m_X & = & \max \left\{ (c d)^{d^2}, (d c / 2)^{d^2 / 2}, (d c^2)^{d^2 / 2} \right\}^n \nonumber \\
& = & \max \left\{ c d, (d c / 2)^{1 / 2}, d^{1 / 2} c \right\}^{n d^2} \nonumber \\
& = & \max \left\{ c d, (d c / 2)^{1 / 2} \right\}^{n d^2} \nonumber \\
& = & \left\{ \begin{matrix} (c d)^{n d^2} & \mathrm{if} \quad 2 c d \ge 1 \\ (c d / 2)^{n d^2 / 2} & \mathrm{if} \quad 2 c d < 1.\end{matrix} \right.
\end{eqnarray}

%%%%%%%%%%%%%%%%%%%%

\subsection{Choices for $T_1$, $T_5$, $T_3$ and $T'_3$}

%%%%%%%%%%%%%%%%%%%%

We can take
\begin{equation}
\label{eq:T1def} 
T_1 = \mathrm{max}(1,c)
\end{equation}
by definition of $T_1$ and of $X = X(c)$.

By definition of $m_X$ we know that
\[
|\mathrm{Norm}_v(y_v)|_v \le m_X^{[k_v : \mathbb{R}] / n}
\]
for every $v \in V_\infty$ and $y = \prod_{v \in V_\infty} y_v \in X$. Here $\mathrm{det}_v(y_v)^d = \mathrm{Norm}_v(y_v)$ so
\[
|\mathrm{det}_v(y_v)|_v^{1 / [k_v : \mathbb{R}]} \le m_X^{1 / (d n)}
\]
for $v$ and $y$ as above. It follows from the definition of $T_5$ that we can take
\begin{equation}
\label{eq:T5def}
T_5 = \mathrm{max}(1,m_X^{1 /(d n)}).
\end{equation}

We chose $\mathcal{D}$ such that $\mathcal{D}_v$ is $\mathrm{M}_{m(v)}(U_v)$ for every $v \in S_f$, where $U_v$ is the the unique maximal $O_v$-order in the division algebra $A_v$. It follows that we can take
\begin{equation}
\label{eq:T3def}
T_3 = T'_3 = 1.
\end{equation}

%%%%%%%%%%%%%%%%%%%%

\subsection{Topological generators and the constants $T_2$ and $T_4$}\label{s:topgens}

%%%%%%%%%%%%%%%%%%%%

We now specify a set $P$ of topological generators for $G_S$ which contains the identity element. If $v$ is archimedean, then $B_v^*$ is isomorphic to $\mathrm{GL}_d(k_v)$ or $\mathrm{GL}_{d/2}(\mathbb{H})$. We claim that there is a set $P_\infty$ of topological generators for $G_S \cap (B_{\mathbb{R}}^* \times \{1\})$ consisting of elements of the form $(x, 1) \in (B_{\mathbb{R}}^* \times B_f^*)$ with $x = \prod_{v \in V_\infty} x_v$ and $|x_v|_v = 1$ for all $v \in V_{\infty}$.

Indeed, if $A_v = \mathbb{H}$ or $\mathbb{C}$, then $\mathrm{GL}_{m(v)}(A_v)$ is connected, so any open subset of $B_v$ generates all of $B_v$. Therefore, the only element needed for these places is the identity. If $A_v = \mathbb{R}$, then $B_v$ has two connected components, determined by the sign of the determinant. Here it suffices to take a topological generator in $G_S$ consisting of the matrix $\mathrm{diag}(-1, 1, \dots, 1)$ at $v$ and the identity at all other $v$ (finite or infinite), which suffices since any open set generates the connected component of the identity. Clearly $|x_v|_v = 1$ for every $x \in P_\infty$ and $v \in V$. This proves the claim.

Now consider $v$ in $S_f$. Then $D_v$ is isomorphic to $\mathrm{M}_{m(v)}(U_v)$ and $B_v$ is isomorphic to $\mathrm{M}_{m(v)}(A_v)$. Let $\lambda_v$ be a prime element of $D_v$, so that $\lambda_v U_v = U_v \lambda_v$ is the unique maximal two-sided proper ideal of $U_v$.

Include in $P_v$ the set of elements of the form $(1, \beta) \in (B_{\mathbb{R}}^* \times B_f^*)$ such that $\beta = \prod_{w \in S_f} \beta_w$ has $\beta_w = 1$ unless $w = v$, and $\beta_v$ is either a permutation matrix, an elementary matrix associated to an element of $U_v$, or a diagonal matrix having all diagonal elements equal to $1$ and the remaining diagonal entry in $U_v^*$. As in the proof of Lemma \ref{lem:detbound1}, every element of $B_{\mathbb{R}}^* \times B_f^*$ can be written as the product of an element in the closure of the group generated by $\cup_{v \in S_f} P_v$ times an element $(x, \alpha)$ in which $x \in B_{\mathbb{R}}^*$ and $\alpha = \prod_{v \in V_f} \alpha_v \in B_f^*$ have the property that for each $v \in V_f$, there are integers $z_1, \ldots, z_{m(v)}$ which may depend on $v$ such that $\alpha_v$ is the diagonal matrix with diagonal entries $\lambda_v^{z_1}, \ldots, \lambda_v^{z_{m(v)}}$.

Also, assume that for each $v \in S_f$, $P_v$ contains elements of the form $(x(v), t(v))$, where $t(v) = \prod_{w \in S_f} t(v)_w$ with $t(v)_w = 1$ if $w \ne v$ and $t(v)_v \in B_f^*$ is the diagonal matrix at the place $v$ having one diagonal entry equal to $\lambda_v$ and the others equal to $1$. Let $x(v) \in B_{\mathbb{R}}^*$ be a real scalar $\tau > 0$ times the identity matrix such that
\[
|\mathrm{Norm}_\infty(x(v))| = \tau^{d^2 n} = |\mathrm{Norm}_f(t(v))|^{- 1} = (\#k(v))^d.
\]

Now the $P_w$ for each $w \in S_f$ together with $P_\infty$ determines a set $P$ of topological generators $(x, \beta)$ for $G_S$ such that for each archimedean place $v$ we have 
\begin{equation}\label{eq:xbound}
|x|_v^{d(v) / [k_v : \mathbb{R}]} \le (m'_{S_f})^{1/n},
\end{equation}
where $m'_{S_f} = 1$ if $S_f = \emptyset$ and otherwise 
\begin{equation}\label{eq:Sfprime}
m'_{S_f} = \max\left\{ (\#k(w))^{d(v) / d }: w \in S_f, v \in S_\infty \right\}.
\end{equation}
Define $m_{S_f}$ to be $1$ when $S_f = \emptyset$ and otherwise
\begin{equation}
\label{eq:normbound}
m_{S_f} = \max \{ \mathrm{Norm}(w): w \in S_f \} = \mathrm{max}\{ \# k(w): w \in S_f\}.
\end{equation}
Since $d(v) \le 2$ and $d(v)/d \le 1$ for $v \in V_\infty$ we have
\begin{equation}\label{eq:xbound2prime}
|x|_v^{d(v)/[k_v : \mathbb{R}]} \le (m'_{S_f})^{1 / n} \le m_{S_f}^{q/n} \quad \mathrm{when}\quad 
q =  \mathrm{min}(2 /d, 1).
\end{equation}
Therefore we can choose
\begin{equation}\label{eq:T2def}
T_2 = (m'_{S_f})^{1 / n} \le m_{S_f}^{q / n}.
\end{equation}
Since all of the non-archimedean components of elements of $P$ are integral, we can choose
\begin{equation}\label{eq:T4def}
T_4 = 1.
\end{equation}

%%%%%%%%%%%%%%%%%%%%

\subsection{An upper bound on $T_6$}

%%%%%%%%%%%%%%%%%%%%

Recall that $T_6$ is the maximum over all subsets $W$ of $V_\infty$ of
\[
T_1^{a(W)} \, T_2^{b(W)} \, T_5^{b(V_\infty \smallsetminus W)}
\]
where $a(W) = \sum_{v \in W} [k_v : \mathbb{R}] m(v)$ and $b(W) = \sum_{v \in W} [k_v : \mathbb{R}]$. Since
\[
T_2 = (m'_{S_f})^{1 / n} \ge 1;
\]
\[
b(W) \le b(V_\infty) = n;
\]
\[
T_1 = \max \{1, c\};
\]
\[
T_5 = \max \{ 1, m_X^{1 / d n} \},
\]
we have an upper bound
\begin{eqnarray}\label{eq:est}
& & T_1^{a(W)} \, T_2^{b(W)} \, T_5^{b(V_\infty \smallsetminus W)} \nonumber \\
&\le & \max \{1, c \}^{a(W)} \, m'_{S_f} \, \max \{1, m_X \}^{b(V_\infty \smallsetminus W) / d n}.
\end{eqnarray}
If $2 c d < 1$, then $c < 1$ and $m_X < 1$ by \eqref{eq:mxdef}, so
\[
T_1^{a(W)} \, T_2^{b(W)} \, T_5^{b(V_\infty \smallsetminus W)} \le m'_{S_f} \quad \mathrm{if} \quad 2 c d < 1.
\]

Now suppose that $2 c d \ge 1$, so that $m_X = (c d)^{n d^2}$ by \eqref{eq:mxdef}. Then $m(v) \le d$ for all $v \in V_\infty$, so \eqref{eq:est} gives
\begin{eqnarray}\label{eq:est2}
& & T_1^{a(W)} \, T_2^{b(W)} \, T_5^{b(V_\infty \smallsetminus W)} \nonumber \\
& \le & \max \{1, c \}^{d b(W)} \, m'_{S_f} \, (c d)^{n d^2 b(V_\infty \smallsetminus W) / d n} \nonumber \\
& \le & \min \{1, c \}^{d b(W)} \, c^{d b(W) + d b(V_\infty \smallsetminus W)}\, m'_{S_f} \, d^{d b(V_\infty \smallsetminus W)} \nonumber \\
& \le & c^{n d} \, m'_{S_f} \, d^{d n} \quad \mathrm{if} \quad 2 c d > 1.
\end{eqnarray}
Putting together \eqref{eq:est} and \eqref{eq:est2} gives
\begin{equation}\label{eq:upper6}
T_6 \le m'_{S_f} \, \max \{ 1, (c d)^{n d} \}.
\end{equation}

%%%%%%%%%%%%%%%%%%%%

\subsection{The explicit bound}\label{ssec:explicit}

%%%%%%%%%%%%%%%%%%%%

Collecting all the above choices leads via Theorem \ref{bounded height} to the following result.

%%%%%%%%%%%%%%%%%%%%

\begin{Theorem}\label{bounded height two}
Suppose $B$ is a central simple division algebra of dimension $d^2$ over a number field $k$, $n = [k : \mathbb{Q}]$, and $s$ is the number of real places of $k$ over which $B$ ramifies. Then there is a maximal order $\mathcal{D}$ of $B$ and functions $f_1(n, d)$ and $f_2(n, d)$ of integer variables $n$ and $d$ for which the following is true. Define
\[
e = \frac{2n}{d(2n-s)}.
\]
Then $e \le 1$. Suppose that $S$ is a finite set of places of $k$ containing all the archimedean places and that $S$ contains all finite places $v$ such that
\[
\mathrm{Norm}(v) \le f_1(n, d) \, d_{\mathcal{D}}^e.
\]
Let $m_{S_f}$ be the maximum norm of a finite place in $S$. Then the group $\Gamma_S$ of $S$-units in $B$ with respect to the order $\mathcal{D}$ is generated by the finite set of elements of height bounded above by
\[
f_2(n, d) \, m'_{S_f} \, d_{\mathcal{D}}^e \le f_2(n, d) \, m_{S_f} \, d_{\mathcal{D}}^e,
\]
where $m'_{S_f}$ is as in $\S$\ref{s:topgens}. In particular, for fixed $n$ and $d$, the height bound for the generating set is polynomial in $m_{S_f}$ and $d_{\mathcal{D}}$.
\end{Theorem}

%%%%%%%%%%%%%%%%%%%%

\begin{proof}
It is clear that $e = 1 / d \le 1$ if $s = 0$, so suppose that $s > 0$. Then $s \le n$ and $d \ge 2$ so we again find that $e \le 1$. The rest of the theorem follows immediately from Theorem \ref{bounded height} and the calculations of the previous subsections.
\end{proof}

%%%%%%%%%%%%%%%%%%%%

\begin{Remark}
We now give closed expressions for $f_1(n, d)$ and $f_2(n, d)$ in the case where $c \geq 1$, with $c$ as in \eqref{eq:cdef}. We leave the adjustments when $c < 1$ as an exercise. We have
\begin{equation}\label{eq:f1def}
f_1(n, d) = d^{n d} \left( \frac{2}{\pi} \right)^{\frac{n d r_2}{n - s / 2}} \left( \frac{2 \sqrt{2}}{\pi} \right)^{\frac{n d s}{2 n - s}}
\end{equation}
\begin{equation}\label{eq:f2def}
f_2(n, d) = d^{n d + n + s} \left( (d - 1)! \right)^{n - s} 2^{s \frac{d^2 - 2 d - 4}{2}} \left( \frac{2}{\pi} \right)^{\frac{n d r_2}{n - s / 2}} \left( \frac{2 \sqrt{2}}{\pi} \right)^{\frac{n d s}{2 n - s}}.
\end{equation}
\end{Remark}

%%%%%%%%%%%%%%%%%%%%

\begin{Remark}
Suppose that $B = k$. If $c < 1$, then $\mathrm{Vol}(X(c)) = 2^n d_\mathcal{D}$ and Minkowski's theorem would imply that $\mathcal{D} = O_k$ contains a non-zero element with norm to $\mathbb{Z}$ less than $1$ in absolute value, which is impossible. Hence $c \ge 1$ and Theorem \ref{bounded height two} exactly reproduces Lenstra's result.
\end{Remark}

%%%%%%%%%%%%%%%%%%%%

\section{An explicit example}\label{sec:Hamilton}

%%%%%%%%%%%%%%%%%%%%

In this section, we compute an explicit example of the bounds in Theorem \ref{bounded height}. Let $B$ denote the quaternion algebra over $\mathbb Q$ ramified exactly at $\{\infty, 2\}$, and let $\mathcal D$ be the Hurwitz order
\[
\mathbb Z \left[ I, J, \frac{1 + I + J + I J}{2} \right],
\]
where $I^2 = J^2 = -1$. Since $\mathbb{R} \otimes_{\mathbb{Q}} B \cong \mathbb{H}$, we have that the Tamagawa measure on $B$ is $4 dx_1dx_2dx_3dx_4$ with respect to the basis $\{1, I, J, I J\}$, and $d_{\mathcal{D}} = 2$.

Recall that
\[
X(c) = \{ x \in \mathbb{H}\ :\ |x|_\infty^{d(v) / [k_v : \mathbb{R}]} \leq c \} = \{ x \in \mathbb{H}\ :\ |x|_\infty^2 \leq c \},
\]
where $|x|_\infty = |\mathrm{N}_\infty(x)|^{1 / d(v)}$. Since $d(v) = 2$ and
\[
\mathrm{N}_\infty(a + b I + c J + d I J) = a^2 + b^2 + c^2 + d^2,
\]
we see that
\[
X(c) = \{ x = a + b I + c J + d I J \in \mathbb{H}\ :\ \|x\| \leq \sqrt{c} \},
\]
where $\|\ \|$ is the usual norm on $\mathbb{R}^4$ with respect to the basis $\{1, I, J, I J\}$. In other words, $X(c)$ is the ball of radius $\sqrt{c}$.

We then see that $X(c)$ has volume
\[
4 \frac{\pi^2 \sqrt{c}^4}{\Gamma(3)} = 2 \pi^2 c^2
\]
with respect to the Tamagawa measure on $\mathbb{H}$, where $\Gamma(s)$ is the usual Gamma function. Since we want
\[
\mathrm{Vol}(X(c)) \geq 2^{\dim_{\mathbb{Q}}(B)} d_{\mathcal{D}} = 32,
\]
we take $c = 4 / \pi$. The constant $m_X$ is the largest square-norm of an element of $X(c)$, which is $16 / \pi^2$.

%%%%%%%%%%%%%%%%%%%%

\subsection{$S = \{\infty\}$}

%%%%%%%%%%%%%%%%%%%%

Since $\sqrt{m_X} < 2$, we see that Theorem \ref{bounded height} applies to \emph{any} set $S$ containing $\{\infty\}$. That is, the set of finite places that must be in $S$ is empty. In the case $S = \{\infty\}$, we use the above to see that
\[
T_2 = T'_3 = T_3 = T_4 = 1, \quad T_1 = T_5 = T_6 = 4 / \pi.
\]
Since $n = 1$ and $d = 2$, plugging these into Theorem \ref{bounded height} gives a height bound of $\frac{4}{\pi} < 2$. The height of an element of $\mathcal{D}$ is its reduced norm, so the elements of height less than $2$ are those with reduced norm $1$. These are the elements of the unit group
\[
\mathcal{D}^* = \left\{ \pm 1, \pm I, \pm J, \pm I J, \frac{\mp 1 \mp I \mp J \mp I J}{2} \right\},
\]
which is well-known to be the binary tetrahedral group (see Theorem 3.7 of \cite{Vigneras}).

%%%%%%%%%%%%%%%%%%%%

\subsection{$S = \{\infty\} \cup\{\ell_i\}_{i = 1}^h$ for a finite set $\{\ell_i\}_{i = 1}^h$ of odd primes}

%%%%%%%%%%%%%%%%%%%%

In this section, we prove Theorem \ref{thm:HamiltonIntro}, the statement of which we now recall.

%%%%%%%%%%%%%%%%%%%%

\begin{Theorem}\label{thm:Hamilton}
Let $B$ be Hamilton's quaternion algebra over $\mathbb{Q}$, that is, the rational quaternion algebra with basis $\{1, I, J, IJ\}$ such that $I^2 = J^2 = 1$ and $I J = - J I$. Let $\mathcal{D}$ be the maximal order
\[
\mathbb{Z} \left[ 1, I, J, \frac{1 + I + J + I J}{2} \right]
\]
and $S = \{\infty, \ell_1, \dots, \ell_h\}$ be a set of places containing the archimedean place $\infty$ and any set $\{\ell_i\}_{i = 1}^h$ of distinct odd primes. Then the unit group $\mathcal{D}_S^*$ is generated by the finite set of elements with reduced norm in $\{1, \ell_1, \dots, \ell_h\}$.
\end{Theorem}

%%%%%%%%%%%%%%%%%%%%

\begin{proof}
Let $\ell_h$ be the largest element of $\{\ell_i\}_{i = 1}^h$. Since all the $\ell_i$ are unramified in $B$ we see that $m_{S_f} = m'_{S_f} = \ell_h$. We can take $T_2 = \ell_h$ and $T_4 =1$ by \eqref{eq:T2def} and \eqref{eq:T4def}, and we can also take $T_3 = T'_3 = 1$ and $T_1 = T_5 = 4/\pi$. This gives $T_6 = \ell_h 4/\pi$, and we conclude that $\mathcal{D}^*$ is generated by elements of height bounded by 
\[
\frac{4}{\pi} \, \ell_h.
\]

We can be more explicit by going back to the statement of Lemma \ref{lem:top-gens}. Recall that
\begin{equation}\label{eq:fund-domain2}
F_X = \{ (x, \beta) \in G_S\ :\ x \in X,~\beta \mathcal{D} \subseteq \mathcal{D},~[\mathcal{D} : \beta \mathcal{D}] \leq m_X \}.
\end{equation}
Since $m_X = 16 / \pi^2 < 2$, we see that if $(x, \beta) \in F_X$ then $\beta \mathcal{D} = \mathcal{D}$. Thus if $\ell$ is an element of the set $P$ of topological generators specified in \S \ref{s:topgens}, and $\gamma \in \Gamma_S \cap (F_X P F_X^{-1})$ as in Lemma \ref{lem:top-gens}, then for each finite place $w$ there are units $u_w, u'_w \in D_w^*$ such that $\gamma_w = u_w \ell_w (u'_w)^{-1}$. Thus $N_w(\gamma)$ equals a unit in $O_w^*$ times $N_w(\ell_w)$ when $N_w:B_w \to k_w$ is the reduced norm. By definition of $P$, this means that the reduced norm $N(\gamma)$ lies in $\{1, \ell_1, \ldots, \ell_h\}$. Therefore $\Gamma_S = \mathcal{D}_S^*$ is generated by set of $\gamma \in \mathcal{D}$ such that $N(\gamma) \in \{1, \ell_1, \ldots, \ell_h\}$.
\end{proof}

The fact that the set $F_X$ defined in the proof of Theorem \ref{thm:Hamilton} is a fundamental domain for the action of $\mathcal{D}_S^*$ on $G_S$ implies the following.

\begin{Corollary}\label{cor:Hamilton}
Let $B$, $\mathcal{D}$, and $S$ be as in Theorem \ref{thm:Hamilton}. For each $1 \le i \le h$, let $T_i$ be the Bruhat--Tits tree associated with $\mathrm{GL}_2(\mathbb{Q}_{\ell_i})$, and set
\[
T = \prod_{i = 1}^h T_i.
\]
Then the action of $\mathcal{D}_S^*$ on $T$ is vertex transitive. That is, given any pair $(x_1, \dots, x_h)$ and $(y_1, \dots, y_h)$, where each $x_i$ and $y_i$ is a vertex of the tree $T_i$, there exists an element $\gamma \in \mathcal{D}_S^*$ so that $\gamma(x_i) = y_i$ for every $1 \le i \le h$.
\end{Corollary}

In comparison, we note that Mohammadi and Salehi Golsefidy \cite{MSG} classified the maximal discrete subgroups with vertex transitive action on the Bruhat--Tits building $\mathfrak{B}$ of a simply connected absolutely almost simple $k$-group over a nonarchimedean field $k$ when $\mathfrak{B}$ has dimension at least $4$. In particular, there are finitely many such maximal groups.

Members of the 2012 Arizona winter school on arithmetic geometry used Corollary \ref{cor:Hamilton} to produce presentations for the groups $\mathcal{D}_S^*$ for various $S$. See \cite{Bourdon}. Unpublished work of Fritz Grunewald computed presentations of $\mathcal{D}_S^1$ for some small $S$, but to our knowledge, \cite{Bourdon} gives the first presentations of groups acting faithfully, irreducibly, and cocompactly on a product of Bruhat--Tits buildings. One such presentation is the following.

%%%%%%%%%%%%%%%%%%%%

\begin{Theorem}[\cite{Bourdon}]\label{thm:3-5}
With $B$ and $\mathcal{D}$ as above, the group $\mathcal{D}_{\{\infty, 3, 5\}}^*$ has presentation with generators (the images in $\mathcal{D}_{\{\infty, 3, 5\}}^*$ of):
\begin{align*}
a &= -1 + I - J - 3 I J \\
b &= -9 - 7 I - J + 7 I J
\end{align*}
and relations:
\begin{eqnarray}
r_1 &=& (b^{-1} a^{-1} b a^{-1})^3 \nonumber \\
r_2 &=& (b^{-1} a^{-2} b a^{-1} b^{-1} a^{-1})^2 \nonumber \\
r_3 &=& (a^{-1} b^{-1} a^{-1} b^{-1} a^{-1} b a^{-1})^2 \nonumber \\
r_4 &=& b^{-1} a b a b^{-1} a^{-1} b^2 a b^{-1} a b a^2 b^{-1} a b a b^{-1} a^2 b a^2 b^{-1} a^{-1} b a^{-2} b^{-1} a^{-2} \nonumber \\
r_5 &=& (b a^2 b^{-1} a b a^{-1} b)^2 \nonumber \\
r_6 &=& b^{-1} a^3 b a^2 b^{-1} a b^{-1} a^{-2} b a^{-1} b^{-1} a \nonumber \\
r_7 &=& b^{-2} a^{-1} b a^{-1} b^{-1} a b a^2 b^{-2} a^{-2} b a^{-1} \nonumber \\
r_8 &=& a b^{-1} a^2 b a^{-1} b^{-1} a^{-2} b a^{-2} b^{-1} a b a \nonumber
\end{eqnarray}
\end{Theorem}

%%%%%%%%%%%%%%%%%%%%

\bibliographystyle{plain}

\bibliography{Lenstra}

%%%%%%%%%%%%%%%%%%%%

\end{document}